\documentclass[10pt]{siamltex}
\usepackage{amsfonts}
\usepackage{latexsym,epic,eepic}
\usepackage[dvips]{graphics,epsfig,color}
\usepackage[color]{showkeys}

\newcounter{cst}
\def \ctel#1{{C_{\refstepcounter{cst}\label{#1}\thecst}\, }}
\def \cter#1{{C_{\ref{#1}}\, }}

\newcounter{cexp}
\def \terml#1{T_{\refstepcounter{cexp}\label{#1}\thecexp}}
\def \termr#1{T_{\ref{#1}}}

\newtheorem{remark}{Remark}[section]

\def\dsp{\displaystyle}

\def\be{\begin{equation}}
\def\ee{\end{equation}}

\def\ba{\begin{array}{lllll}}
\def\ea{\end{array}}

\def\beqsys {\be\ba \left \{ \begin{array}{l}}
\def\eeqsys {\end{array} \right . \ea\ee }

\def\beqsysno {\be\ba \left \{ \begin{array}{l}}
\def\eeqsysno {\end{array} \right . \ea\ee}

\def\centers{{\cal P}}

\def\cv{K}

\def\d{{\rm d}}

\def\disc{{\cal D}}

\def\dcvedge{d_{\cv,\edge}}

\def\diam{\hbox{\rm diam}}
\def\div{{\rm div}}

\def\dkl{d_{\kl}}

\def\dr{\partial}
\def\dt{\delta\!t}

\def\edge{\sigma}
\def\edges{{\cal E}}
\def\edgesint{{\cal E}_{{\rm int}}}
\def\edgesext{{\cal E}_{{\rm ext}}}
\def\edgecvcvv{K|L}
\def\edgescv{{\cal E}_K}

\def\eps{\varepsilon}

\def\grad{\nabla}

\def\half{{\frac 1 2}}

\def\kl{{K|L}}

\def\lap{\Delta}

\def\mcv{\meas_K}
\def\meas{{\rm m}}
\def\medge{\meas_{\edge}}
\def\mesh{{\cal M}}

\def\mK{\meas_K}
\def\mkl{\meas_{K\vert L}}

\def\n{\mathbf{n}}

\def\N{\mathbb{N}}
\def\NN{{\cal N}}

\def\normedeu(#1){\|#1\|_{L^2(\Omega)}}

\def\O{\Omega}

\def\phi{\varphi}

\def\R{\mathbb{R}}
\def\refe#1{(\ref{#1})}

\def\regul{\hbox{\rm regul}}

\def\s{\sigma}

\def\size{\hbox{\rm size}}

\def\sumi{\dsp{\sum_{i=1}^{d}}}
\def\sumj{\dsp{\sum_{j=1}^{d}}}
\def\sumk{\dsp{\sum_{k=1}^{d}}}

\def\tends{\to}
\def\tkl{\frac {\mkl}{\dkl} }
\def\tks{\frac {\medge}{\dcvedge}}

\title{Convergence analysis of a colocated  finite volume scheme
 for the 
 incompressible Navier-Stokes equations on 
general 2 or 3D meshes }

\author{
R. Eymard \thanks{Universit\'e de Marne-la-Vall\'ee, France,
({\tt eymard@math.univ-mlv.fr})}
\and
R. Herbin \thanks{Universit\'e de Provence, France
({\tt herbin@cmi.univ-mrs.fr})}
\and J.C. Latch\'e \thanks{DPAM, Institut de Radioprotection et 
Suret\'e Nucl\'eaire,
({\tt jean-claude.latche@irsn.fr})}
}

\begin{document}

\maketitle

\begin{abstract}
We study a colocated cell centered finite volume method for the approximation 
of the incompressible Navier-Stokes  equations posed on a 2D or 3D finite domain.
The discrete unknowns are the components of the velocity and the pressures, all of them 
colocated at the center of the cells of a unique mesh; hence the need for a stabilization 
technique,   which we choose of the Brezzi-Pitk\"aranta type.
The scheme features two essential properties: the discrete gradient is the transposed of 
the divergence terms and the discrete trilinear form associated to nonlinear advective terms 
vanishes on discrete divergence free velocity fields.
As a consequence, the scheme is proved to be unconditionally stable and convergent for the 
Stokes problem, the steady and the transient Navier-Stokes equations.
In this latter case, for a given sequence of approximate solutions computed on meshes the size of which 
tends to zero, we prove, up to a subsequence, the $L^2$-convergence of the components of the velocity, 
and, in the steady case, the weak $L^2$-convergence of the pressure.
The proof relies on the study of space and time translates of approximate solutions, which
allows the application of Kolmogorov's theorem.
The limit of this subsequence is then shown to be a weak solution
of the Navier-Stokes equations.
Numerical examples are performed to obtain numerical convergence rates in both the linear
and the nonlinear case.
\end{abstract}



\section{Introduction}

We are interested in this paper in finding an approximation of the  
 fields $\bar u = (\bar u^{(i)})_{i=1,\ldots,d}$~:~$\Omega \times [0,T]
 \rightarrow \R^d$, and $\bar p$~:~$\Omega \times [0,T]
 \rightarrow \R$,
  weak solution to the incompressible Navier-Stokes equations which write:
\be
\qquad \ba\dsp
\partial_t \bar u^{(i)} -\nu \lap \bar u^{(i)} + \partial_i \bar p + \sumj \bar u^{(j)}
\partial_j \bar u^{(i)}
 = f^{(i)}
\hbox{ in } \Omega\times(0,T), \ \mbox{ for } i = 1,\ldots,d,\\ \dsp
\div \bar u = \sumi \partial_i \bar u^{(i)}   = 0 \hbox{ in } \Omega\times(0,T).
\ea\label{nstocontt}
\ee
with a homogeneous Dirichlet boundary condition for $\bar u$ and the initial condition
\be
\ba\dsp
\bar u^{(i)}(\cdot,0) = \bar u^{(i)}_{\rm ini} \hbox{ in } \Omega \mbox{ for } i =
1,\ldots,d.
\ea\label{nstoconti}
\ee
In the above equations,  $\bar u^{(i)}$, $i=1,\ldots,d$ denote the components of the velocity of a fluid which flows in a domain $\O$ during the time $(0,T)$, $\bar p$ denotes the
pressure, $\nu>0$ stands for the viscosity of the fluid.
We make the following assumptions:
\begin{eqnarray}
&\qquad &
\O \mbox{ is a polygonal open bounded connected subset of }\R^d,\ d=2
\mbox{ or }3,
\label{hypomegat}
\\&&
T > 0\mbox{ is the finite  duration of the flow,}
\label{hyptimet}
\\&&
\nu \in (0,+\infty),
\label{hypnut}
\\&&
\bar u_{\rm ini}\in L^2(\O)^d,
\label{hypuini}
\\&&
f^{(i)} \in L^2(\O\times (0,T)), \ \mbox{ for } i = 1,\ldots,d.
\label{hypfgt}
\end{eqnarray}

We denote by $x = (x^{(i)})_{i=1,\ldots,d}$ any point of $\O$, by $\vert .\vert$ the Euclidean norm in $\R^d$, {\it i.e.}:
$
\vert x \vert^2 = \sumi (x^{(i)})^2
$
and by $\d x$ the $d$-dimensional Lebesgue measure $\d x = \d x^{(1)}\ldots \d x^{(d)}$.

\medskip
The weak sense that we consider for the Navier-Stokes equations is the following.

\medskip
\begin{definition}[Weak solution for the transient Navier-Stokes equations]\label{weaksolt}\\
Under hypotheses \refe{hypomegat}-\refe{hypfgt}, let the function space $E(\O)$ be defined by:
\begin{eqnarray}
E(\O) := \{ \bar v=(\bar v^{(i)})_{i=1,\ldots,d} 
\in H^1_0(\O)^d, \div \bar v  = 0
\mbox{ a.e. in
}\O\}.\label{e0}
%,
%\qquad \nonumber
%\\
%\mbox{\it i.e. } \int_O (\div \bar v)(x)\ q(x)\, \d x=0 \quad \forall q \in L^2(0) \}.
\end{eqnarray}
Then  $\bar u$ is called a weak solution of \refe{nstocontt}-\refe{nstoconti} if
$\bar u\in L^2(0,T;E(\O))\cap L^\infty(0,T;L^2(\O)^d)$ and:
\be
\left\{\ba \dsp
\forall \varphi \in L^2(0,T;E(\O)) \cap C^\infty_c(\O\times (-\infty,T))^d,
\\[3ex] \dsp \qquad
-\int_0^T\int_\O \bar u(x,t) \cdot\partial_t\varphi(x,t)\, \d x\, \d t
- \int_\O \bar u_{\rm ini}(x) \cdot\varphi(x,0)\, \d x
\\ \dsp \qquad
+ \nu \int_0^T\int_\O \nabla \bar u(x,t):\nabla \varphi(x,t)\, \d x\, \d t \
+ \int_0^T b(\bar u(\cdot,t),\bar u(\cdot,t),\varphi(\cdot,t))\, \d t
\\ \dsp \hfill
= \int_0^T\int_\O f(x)  \cdot\varphi(x,t)\, \d x\, \d t
\ea\right.
\label{nstocontft}
\ee
where, for all $\bar u,\bar v\in H^1_0(\O)^d$ and for a.e. $x\in\O$, we use the following notation:
\[
\nabla \bar u(x):\nabla \bar v(x) = \sumi \nabla \bar u^{(i)}(x)\cdot\nabla \bar v^{(i)}(x)
\]
and where the trilinear form $b(.,.,.)$ is defined, for all $\bar u,\bar v,\bar w \in (H^1_0(\O))^d$, by
\be
b(\bar u,\bar v,\bar w) = \sumk \sumi \int_\O \bar u^{(i)}(x) \partial_i \bar v^{(k)}(x) \bar w^{(k)}(x) \, \d x.
\label{deftricont}
\ee
\end{definition}

\medskip
\begin{remark}
From \refe{nstocontft}, we get that a weak solution $u$ of 
\refe{nstocontt}-\refe{nstoconti} in the sense of Definition \ref{weaksolt} satisfies $\partial_t \bar u\in L^{4/d}(0,T;E(\O)')$, and is therefore a weak solution in the classical sense,
such that $\bar u(\cdot,0)$ is the orthogonal $L^2$-projection of $\bar u_{\rm ini}$ 
on $\{\bar v\in L^2(\O)^d, \div \bar v = 0, {\rm trace}(\bar v\cdot n_{\dr\O},\dr\O) = 0\}$
(see for example \cite{temam} or \cite{bf}).
\end{remark}

\bigskip

Numerical schemes for the Stokes equations  and the 
Navier-Stokes equations   have been extensively studied:
see
 \cite{giraultraviart, patankar,peyret-taylor,pironneau, gunzburger,glo} and references therein.
Among different schemes, finite element schemes and finite volume schemes
are frequently used for mathematical or engineering studies.
An advantage of  finite volume schemes is that the unknowns
are approximated by piecewise constant functions: this makes it  easy to take into account
additional nonlinear phenomena or the coupling with algebraic or differential equations,
 for instance in the case
of reactive flows; in particular, one can find in 
\cite{patankar} the presentation of the
classical finite volume scheme on rectangular meshes, which has been the basis
of many industrial applications.
However, the use of rectangular grids
makes an important limitation to the type of domain which can be gridded
and more recently, finite volume schemes for the Navier-Stokes equations on 
triangular grids have been presented:
see for example \cite{nico} where the vorticity formulation is used and 
\cite{boivin} where primal variables 
are used with a Chorin type projection method to ensure the divergence condition. 
Proofs of convergence for finite volume type schemes for the Stokes and steady-state
Navier-Stokes  equations  are have recently been given 
for staggered grids
\cite{chou}, \cite{nico}, \cite{cras}, \cite{stagg}, \cite{benartzi}, 
following the pioneering work  of Nicolaides {\it et al.} \cite{nicol}, \cite{nicwu}. 

\medskip
In this paper, we propose the mathematical and numerical analysis
of a discretization method which uses the primitive variables, that is the velocity 
and the pressure, both 
approximated by piecewise constant functions 
on the cells of a 2D or 3D mesh. We emphasize that the approximate velocities and 
pressures are colocated, and therefore, no dual grid is needed. 
The only requirement on the mesh  is a geometrical 
assumption needed for the consistency of the approximate diffusion flux 
(see \cite{book} and
section \refe{secdisc} for a precise definition of the admissible discretizations).

\medskip

As far as we know, this work is a first proof of the convergence,  of 
a finite volume scheme which is of large
interest in industry.  Indeed, industrial CFD codes (see e.g. \cite{fluent}, \cite{neptune}) use 
colocated cell centered 
finite volume schemes; leaving aside implementation considerations, the principle of these
schemes  seems
to differ from the present scheme only by the stabilization choice. The main reasons why this scheme is
so popular in industry are:
\begin{itemize}
\item a colocated arrangement of the unknowns,
\item a very cheap assembling step,
(no numerical integration to perform)
\item an easy coupling with other systems of equations.
\end{itemize}

\medskip

The finite volume scheme studied here is based on three basic ingredients.
First, a stabilization technique {\it \`a la } Brezzi-Pik\"aranta \cite{brez} is used
 to cope with the instability of colocated velocity/pressure approximation spaces.
Second, the discretization of the pressure gradient in the momentum balance equation 
is performed to ensure, 
by construction, that it is the transpose of the divergence term of the continuity constraint.
Finally, the contribution of the discrete nonlinear advection term 
to the kinetic energy balance vanishes for discrete divergence free velocity fields, 
as in the continuous case. 
These features appear to be essential in the proof of convergence.

\medskip
%As performed for a wide range of problems in \cite{book}, we  present here the finite volume scheme at
% hand in a discrete functional framework, defining the discrete analogue of the continuous Laplace, gradient, divergence and transport operators, each of them enjoying properties similar to their continous counterparts.
%In return, we are able to recast our finite volume discretization in a variational-like framework, which reads for the transient Navier-Stokes equations, to give an example:
%\[
%\ba \dsp
%\mbox{Find } u_{n+1}\in E_\disc(\O) \mbox{ such that, } \forall  v \in E_\disc(\O),\ \forall n \in \N,
%\\[2ex] \dsp
%\int_\O (u_{n+1}(x)-u_{n}(x))\cdot v(x) \d x + \nu  \dt [u_{n+\half},v]_{\disc}
%\\ \dsp \hspace{20ex}
%+ \dt\ b_\disc(u_{n+\half},u_{n+\half},v)= \int_{n\dt}^{(n+1)\dt}\int_\O f(x,t)\cdot v(x)\,\d x\, \d t
%\ea
%\]
%where $u_{n+\half} = \half(u_{n+1}+u_n)$ and $ E_\disc(\O)$ is the space of discretely divergence free functions, endowed with the inner product $ [\cdot,\cdot]_{\disc}$ (see definition \ref{espdis}).

\medskip
We are then able to prove the stability of the scheme and the convergence of discrete solutions 
towards
 a solution of the continuous problem when the size of the mesh tends to zero, 
 for the steady linear 
 case (generalized Stokes problem), the stationary and the transient Navier-Stokes equations, 
 in 2D and 3D.
Our results are valid for general meshes, do  not require any assumption on 
the regularity of 
the continuous solution nor, in the nonlinear case, any small data condition.
We emphasize that the convergence of the 
fully discrete (time and space) approximation is proven here, 
using an original estimate on the time 
translates, which yields, combined with a classical estimate 
on the space translates, 
a sufficient relative compactness property.

\medskip
An error analysis is performed in the steady linear case, under regularity assumptions on 
the solution.
An error bound of order 0.5 with respect to the step size is obtained  in the discrete 
$H^1$ norm and 
the $L^2$ norm for  respectively the velocity and the pressure. Of course, this is probably 
not a sharp estimate, as can
be seen from the numerical results shown in Section \ref{secnum}.
Indeed, a better rate of convergence can be proved under  
additional assumptions on the mesh  \cite{EHL}.
 
\medskip
%Our numerical experiments show that the scheme studied here shares qualitatively the 
%same convergence properties as by now standard discretisations, as the so-called $P1/P1$ finite element 
%scheme or, in the finite volume context, the celebrated MAC scheme \cite{harlow}, with some additional 
%desirable features: to mention some, a colocated arrangement of the unknowns 
%(which is used in most industrial codes, see \cite{fluent}, \cite{neptune}), a very cheap assembling step 
%(no numerical integration to perform), especially for problems with complex (for instance, 
%tabulated) constitutive laws and, finally, the possibility to couple, within the finite volume context, 
%this Navier-Stokes solver with by-construction monotone discretizations of transport equations 
%of passive scalars.

\bigskip
This paper is organized as follows.
In section \ref{secdisc}, we introduce the discretization tools together with some  
 discrete functional analysis tools. 
Section \ref{secfvslin} is devoted to the linear steady problem (Stokes problem), 
for which the finite volume scheme is given and  convergence analysis and error estimates are detailed.
The complete finite volume scheme for the nonlinear case is presented in section \ref{secfvsnlin}, in both
the steady and transient cases.
We then develop the analysis of its convergence to a weak solution of the continuous problem.
We  give some numerical results in section \ref{secnum}, and finally conclude with some remarks 
on open problems (section \ref{seconcrem}).

\section{Spatial discretization and discrete functional analysis}\label{secdisc}

\subsection{Admissible discretization of $\O$}

We first recall the notion of admissible discretization for a finite volume method, which
is given in \cite{book}.

\medskip
\begin{definition}[Admissible discretization, steady case]\label{adisc}
Let $\O$ be an open bounded  polygonal (polyhedral if $d=3$) subset of $\R^d$, and $\dr \O =  \overline{\O}\setminus\O$ its boundary.
An admissible finite volume discretization of $\O$, denoted by $\disc$, is given by $\disc=(\mesh,\edges,\centers)$, where:
\begin{itemize}

\item[-] $\mesh$ is a finite family of non empty open polygonal convex disjoint subsets of
$\O$ (the ``control volumes'') such that $\overline{\O}= \dsp{\cup_{K \in \mesh} \overline{K}}$.
For any $K\in\mesh$, let $\dr K  = \overline{K}\setminus K$ be the boundary of $K$
and $\mK>0$ denote the area of $K$.

\item[-] $\edges$ is a finite family of disjoint subsets of $\overline{\O}$ (the ``edges'' of the mesh), such that, for all $\sigma\in\edges$, there exists a hyperplane $E$ of $\R^d$ and $K\in\mesh$ with $\overline{\sigma} = \dr K \cap E$ and $\sigma$ is a non empty open subset of $E$.
We then denote by $\medge>0$ the (d-1)-dimensional measure of $\sigma$.
We assume that,for  all $K \in \mesh$, there exists  a subset $\edgescv$ of $\edges$
such that $\dr K  = \dsp{\cup_{\sigma \in \edgescv}}\overline{\sigma} $.
It then results from the previous hypotheses that, for all $\sigma\in\edges$, either $\sigma\subset \dr\O$ or there exists $(K,L)\in \mesh^2$ with $K \neq L$ such that
$\overline{K} \cap \overline{L} = \overline{\sigma}$;
we denote in the latter case $\sigma = K|L$.

\item[-] $\centers$ is a family of points of $\O$ indexed by $\mesh$, denoted by $\centers = (x_K)_{K \in \mesh}$.
The coordinates of $x_K$ are denoted by $x^{(i)}_K$, $ i = 1,\ldots,d.$
The family $\centers$ is  such that, for all $K \in \mesh$, $x_K \in K$. Furthermore, for all $\sigma\in\edges$ such that there exists  $(K,L)\in \mesh^2$ with $\sigma = K|L$, it is assumed  that  the straight line $(x_K,x_L)$ going through $x_K$ and $x_L$ is orthogonal to $K|L$.
For all $K \in \mesh$ and all $\sigma \in \edgescv$, let $z_\sigma$ be the orthogonal projection of $x_K$ on $\sigma$.
We suppose that $z_\sigma\in\sigma$.
\end{itemize}
 
\end{definition}

\begin{figure}[ht]
\begin{center}
%DESSIN 3
\input{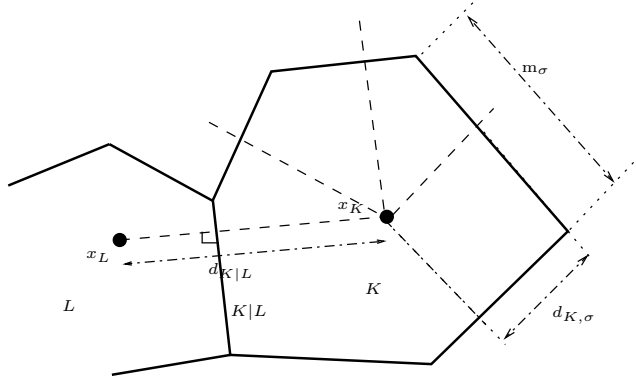}
\end{center}
\caption{Notations for an admissible mesh}\label{fig_maille}
\end{figure}

An example of two neighbouring control volumes $K$ and $L$ 
 of $\mesh$ is depicted in
 Figure \ref{fig_maille}.  

The following notations are used.
The size of the discretization is defined by:
\[
\size(\disc)= \sup\{\hbox{\rm diam}(K), K\in \mesh\}.
\]
For all $K \in \mesh$ and $\sigma \in \edgescv$, we denote by
 $\n_{K,\sigma}$ the unit vector normal to $\sigma$ outward to $K$.
We denote by $d_{K,\sigma}$ the Euclidean distance between $x_K$ and $\sigma$.
The set of interior (resp. boundary) edges is denoted by $\edgesint$ 
(resp. $\edgesext$), that is $\edgesint = \{\sigma \in \edges;$ $\sigma \not \subset \partial \O \}$ (resp.  $\edgesext = \{\sigma \in \edges;$ $\sigma \subset \partial \O \}$).
For all $K\in\mesh$, we denote by $\NN_K$ the subset of $\mesh$ of the 
neighbouring control volumes.
For all $K\in\mesh$ and $L\in\NN_K$, we set $\n_{KL} =\n_{K,K|L}$, 
we denote by $\dkl$
the Euclidean distance between $x_K$ and $x_L$.

We shall measure the regularity of the mesh through the function $\regul(\disc)$ defined by
\be
\ba \regul(\disc) = \inf&\left\{\frac {d_{K,\sigma}} {{\rm diam} (K)}, \
 K\in \mesh,\ \sigma\in\edgescv \right\}\\ 
 & \cup \left\{\frac {d_{K,K|L}} {\dkl}, \
 K\in \mesh,\ L \in \NN_K \right\}\cup
\left\{ \frac 1 {{\rm card}(\edgescv)},\  K\in \mesh \right\}.
\label{regul}
\ea \ee

\subsection{Discrete functional properties}

Finite volume schemes are  discrete balance equations with an adequate 
approximation of the
fluxes, see e.g.  \cite{book}. Recent works dealing with cell centered finite volume methods 
for
elliptic problems \cite{convpardeg}, 
\cite{cras100}, 
\cite{stagg} introduce an equivalent variational formulation in adequate 
functional spaces. 
Here we shall follow this latter path, also introducing discrete analogues 
of the continuous Laplace, 
gradient, divergence and transport operators, 
each of them featuring properties similar to their continuous counterparts.

\begin{definition}\label{espdis}
Let $\O$ be an open bounded  polygonal subset of $\R^d$, with $d\in \N_*$.
Let $\disc=(\mesh,\edges,\centers)$ be an admissible finite volume discretization of $\O$ in the sense of definition \ref{adisc}.
We denote by $H_\disc(\O)\subset L^2(\O)$ the space of functions which are piecewise constant on each control volume $K\in\mesh$.
For all $w\in H_\disc(\O)$ and for all $K\in\mesh$, we denote by $w_K$ the constant value of $w$ in $K$.
The space $H_\disc(\O)$ is embedded with  the following  Euclidean structure:
For $(v,w)\in (H_\disc(\O))^2$, we first define the following inner product 
(corresponding to Neumann boundary conditions)
\be
\langle v,w\rangle_{\disc} = \dsp\frac 1 2 \sum_{K\in\mesh} 
\sum_{L\in\NN_K}
\tkl (v_L - v_K)(w_L - w_K).  
\label{amundis}
\ee
We then define another inner product (corresponding to Dirichlet boundary conditions)
\be
[v,w]_{\disc} = \langle v,w\rangle_{\disc} +\sum_{K\in\mesh} 
\sum_{\sigma\in\edgescv\cap\edgesext}
\tks v_K w_K.
\label{amudis}
\ee

Next, we define a seminorm and a norm in $H_\disc(\O)$ 
(thanks to the discrete Poincar\'e inequality \refe{poindis} given below) by
\[\vert w \vert_{\disc} = \left( \langle w,w\rangle_{\disc} \right)^{1/2}, \qquad
\Vert w \Vert_{\disc} = \left( [w,w]_{\disc} \right)^{1/2}.
\]
We define the interpolation operator $P_\disc~:~C(\O)\to H_\disc(\O)$ by
$(P_\disc\varphi)_K = \varphi(x_K)$, for all $K\in\mesh$, for all $\varphi\in C(\O)$.
%% RE Comme dans vf100, l'operateur d'interpolation est en fait defini dans C(\O)
%% RE j'ai donc change H2 en C ici et un peu plus bas
\medskip

Similarly,  for $u = (u^{(i)})_{i=1,\ldots,d} \in (H_\disc(\O))^d$, $v = (v^{(i)})_{i=1,\ldots,d}\in (H_\disc(\O))^d$ and $w = (w^{(i)})_{i=1,\ldots,d}\in (H_\disc(\O))^d$, we define:
\[ 
\Vert u \Vert_{\disc} = \left( \sumi [u^{(i)}, u^{(i)}]_{\disc} \right)^{1/2},
 \qquad
[v,w]_{\disc} = \sumi[v^{(i)}, w^{(i)}]_{\disc},
\]
and $P_\disc~:~C(\O)^d\to H_\disc(\O)^d$ by $(P_\disc\varphi)_K = \varphi(x_K)$, for all $K\in\mesh$, for all 
$\varphi\in C(\O)^d$.
\end{definition}

The discrete Poincar\'e inequalities (see \cite{book}) write:
\be
\Vert w \Vert_{L^2(\O)} \le \diam(\O) \Vert w \Vert_{\disc}, \ \forall w\in 
H_\disc(\O),
\label{poindis}
\ee
and there exists $C_\O >0$, only depending on $\O$, such that
\be
\Vert w \Vert_{L^2(\O)}^2 \le C_\O \vert w\vert_{\disc}^2, \ \forall w\in 
H_\disc(\O) \hbox{ with }\int_\O w(x)\d x = 0.
\label{poindismoy}
\ee

%% RE ici, il y avait une curieuse duplication
We  define a discrete divergence operator $\div_{\disc}~: (H_\disc(\O))^d \to H_\disc(\O)$,
by:
\be
\div_{\disc}(u)(x) = \frac 1 {\mK}
\sum_{L\in\NN_K}   A_{KL} \cdot (u_K+u_L),\
\mbox{ for a.e. }x\in K,
\forall K\in\mesh,\label{divdisc}
\ee
with
\be
\ba\dsp
A_{KL} =  \tkl\frac {x_L  - x_K } 2 = \half\  \mkl \ \n_{KL},
\ \forall K\in\mesh,\ \forall L\in\NN_K.
\ea\label{defAB}
\ee

We then set $E_\disc(\O) = \{u\in (H_\disc(\O))^d, \div_\disc(u) = 0\}$.

\begin{remark}
Any definition of $A_{KL}$ such that $A_{KL} =  \mkl a_{KL} \n_{KL}$ with $a_{KL} \ge 0$ and $a_{KL}+a_{LK} = 1$, combined with the definition $\div_{\disc}(u)(x) = \frac 1 {\mK}
\sum_{L\in\NN_K}  ( A_{KL} \cdot u_K - A_{LK}  \cdot u_L)$, produces the same results
of convergence as those which are proven in this paper.
On particular meshes, one can prove a better error estimate, choosing $a_{KL} = d(x_L,K|L)/d_{KL}$ (see \cite{EHL}).
Nevertheless, in the general framework of this paper, other choices do not improve the convergence result and the error estimate.
Therefore, we set in this paper $a_{KL} = 1/2$, which corresponds to \refe{defAB}.
The advantage of this choice is that it leads to simpler notations and shorter equations.
\end{remark}

The adjoint of this discrete divergence defines a discrete gradient $\grad_{\disc}~: H_\disc(\O) \to (H_\disc(\O))^d$:
\be 
(\grad_{\disc} u)_K = \frac 1 {\mK}
\sum_{L\in\NN_K}   A_{KL} (u_L - u_K),\
\forall K \in \mesh, \ \forall u\in H_\disc(\O). \label{discgrad}
\ee
This operator $\grad_{\disc}$ then satisfies the following property.
\begin{proposition}\label{weakconvgrad}
Let $(\disc^{(m)})_{m\in\N}$ be a sequence of admissible discretizations of $\O$ in the sense of Definition \ref{adisc}, such that $\dsp\lim_{m\to\infty}\size(\disc^{(m)}) = 0$.
Let us assume that there exists $C >0$ and $\alpha\in [0,2)$ and a sequence $(u^{(m)})_{m\in\N}$ such that $u^{(m)}\in H_{\disc^{(m)}}(\O)$ and 
$ \vert u^{(m)} \vert_{\disc_m}^2 \le C\ \size(\disc^{(m)})^{-\alpha}$, for all $m\in\N$.

Then the following property holds:
\be
\lim_{m\to +\infty} \int_\O  \left(P_{\disc_{m}}\varphi(x) \grad_{\disc_{m}}  u^{(m)}(x) 
 +  u^{(m)}(x)\grad \varphi(x)\right) dx = 0,\ 
  \forall \varphi\in C^\infty_c(\O),
\label{wcvgrad}\ee 
and therefore:
\be
\lim_{m\to +\infty} \int_\O \grad_{\disc_{m}} u^{(m)}(x) \cdot P_{\disc_{m}}\psi(x) \d x = 0,\ 
 \forall \psi\in C^\infty_c(\O)^d\cap E(\O),
\label{wcvdiv}
\ee
where $E(\O)$ is defined by \refe{e0}.
\end{proposition}

\begin{proof} 
Let us assume the hypotheses of the above lemma, and let $i=1,\ldots,d$ and $\varphi\in C^\infty_c(\O)$ be given.
Let us study, for $m\in\N$, the term 
\[  
\terml{wcv1}^{(m)} = 
\int_\O \left(P_{\disc_{m}}\varphi(x) \grad_{\disc_{m}}  u^{(m)}(x) +  
u^{(m)}(x)\grad \varphi(x) \right)
d x.
\]
From \refe{defAB} and \refe{discgrad}, we get that 
\[
\termr{wcv1}^{(m)} = \sum_{\s\in\edgesint,\s = K|L}
(u^{(m)}_L- u^{(m)}_K) \mkl\ R_{KL}^{(m)},  
\]
where
\[
 R_{KL}^{(m)} = 
 \left(\frac 1 2(\varphi (x_K)+\varphi (x_L)) - \frac 1 {\mkl}
\int_{K|L} \varphi (x) \d\gamma(x)\right)\n_{KL}.
\]
Thanks to the Cauchy-Schwarz inequality, 
\[
\vert \termr{wcv1}^{(m)} \vert^2\le \vert u^{(m)} \vert_{\disc_m}^2 
\sum_{\s\in\edgesint,\s = K|L} \left|R_{KL}^{(m)}\right|^2 \mkl d_{KL}.
\]
One has $\sum_{\s\in\edgesint,\s = K|L} \mkl d_{KL} \le d \meas(\O).$
Thanks to the existence of $C_\varphi>0$ which only depends on $\varphi$ such that $|R_{KL}^{(m)}|\le C_\varphi\size(\disc^{(m)})$ and since $\alpha < 2$, we then get that
\[
\lim_{m\to \infty} \termr{wcv1}^{(m)}  = 0,
\]
which yields \refe{wcvgrad}. 
\end{proof}

\begin{proposition}[Discrete Rellich theorem]\label{cpct}
Let $(\disc^{(m)})_{m\in\N}$ be a sequence of admissible discretizations of $\O$ in the sense of definition \ref{adisc}, such that $\dsp\lim_{m\to\infty}\size(\disc^{(m)}) = 0$.
Let us assume that there exists $C >0$ and a sequence  $(u^{(m)})_{m\in\N}$ such that
$u^{(m)}\in H_{\disc^{(m)}}(\O)$ and $ \Vert u^{(m)} \Vert_{\disc_m} \le C$ for all $m\in\N$.

Then, there exists
$\bar u\in H^1_0(\Omega)$ and a subsequence of $(u^{(m)})_{m\in\N}$, again denoted $(u^{(m)})_{m\in\N}$, such that:
\begin{enumerate}
\item the sequence $(u^{(m)})_{m\in\N}$ converges in $L^2(\Omega)$ to $\bar u$ as $m\to +\infty$,
\item for all $\varphi\in C^\infty_c(\O)$, we have
\be
\lim_{m\to +\infty} [u^{(m)},P_{\disc_m}\varphi]_{\disc_m} = \int_\O \grad \bar u(x)\cdot\grad\varphi(x) \d x,
\label{cstd}
\ee
\item $\grad_{\disc_{m}} u^{(m)}$ weakly converges to $\grad \bar u$  in $L^2(\Omega)^d$ as $m\to +\infty$ and \refe{wcvgrad} holds.
\end{enumerate}
\end{proposition}

\begin{proof}
The proof of the first two  items is given in \cite{book} (see proof of Theorem 91. pp
773--774).
Since we have $\vert u^{(m)}\vert_{\disc_m}\le  \Vert u^{(m)} \Vert_{\disc_m}$, we can
apply proposition \ref{weakconvgrad}, which gives the third item.
\end{proof}

\begin{remark}
Following  \cite{hcv}, if we denote
\[
{\cal D}_{K,\sigma} = \{ t x_K + (1-t) y,\ t\in(0,1),\ y\in\sigma\},\ 
\forall K\in\mesh,\ \forall \sigma\in\edgescv,
\]
we may alternatively define a discrete gradient 
$\tilde \grad_{\disc}~: H_\disc(\O) \to (L^2(\O))^d$, by:
\[
\ba \dsp
\hbox{ for all } K\in\mesh,\\
\tilde \grad_{\disc}u(x) = \frac d {d_{KL}} (u_L - u_K)\n_{KL}, \hbox{ for 
a.e. }x\in {\cal D}_{K,K|L}\cup {\cal D}_{L,K|L},
  \ \forall L\in\NN_K, \\ \dsp
 \tilde \grad_{\disc}u(x) = \frac d {d_{K,\sigma}} (0 - u_K)\n_{K,\sigma}, 
\hbox{ for a.e. }x\in {\cal D}_{K,\sigma},
 \ \forall \sigma\in\edgescv\cap\edgesext.
\ea
\]
A result similar to that of Proposition \ref{cpct} holds with this definition of 
a discrete  gradient, and in fact, it can be shown that
the weak convergence of $\tilde \grad_{\disc_m} u^{(m)} $  is equivalent 
to the weak convergence of $  \grad_{\disc_m} u^{(m)} $.
\end{remark}

\section{Approximation of the linear steady problem}\label{secfvslin}

\subsection{The Stokes problem}

We first study the following linear steady problem: find an approximation of  
$\bar u $ and $\bar p$, weak solution to the generalized Stokes equations, which write:
\be
\ba\dsp
\eta \bar u -\nu \lap \bar u + \grad \bar p = f
\hbox{ in } \Omega\\ \dsp
\div \bar u   = 0 \hbox{ in } \Omega,
\ea\label{stocont}
\ee
For this problem, the following assumptions are made:
\be
\O \mbox{ is a polygonal open bounded connected subset of }\R^d,\ d=2 
\mbox{ or }3
\label{hypomega}
\ee
\be
\nu\in (0,+\infty),\ \eta \in [0, + \infty),
\label{hypnu}
\ee
\be
f \in L^2(\O)^d.
\label{hypfg}
\ee
We then consider the following weak sense for problem \refe{stocont}.

\begin{definition}[Weak solution for the steady Stokes equations]\label{weaksol}
 
Under hypotheses \refe{hypomega}-\refe{hypfg}, let $E(\O)$ be defined by \refe{e0}.
Then  $(\bar u,\bar p) $ is called a weak solution of \refe{stocont}  (see e.g. \cite{temam} or \cite{bf}) if
\be
\left\{\ba
\bar u \in E(\O), \ \bar p\in L^2(\O) \hbox{ with } \int_\O \bar p(x) \d x = 0,\\
\dsp \eta \int_\O \bar u(x)\cdot \bar v(x) \d x +
\dsp\nu \int_\O \nabla \bar u(x) :\nabla \bar v(x)
\d x -  \\
\dsp \int_\O \bar p(x)\div \bar v(x) \d x= \dsp \int_\O f(x)\cdot  \bar v(x) \d x,
 \  \forall \bar v \in H^1_0(\O)^d.
\ea\right.\label{stocontf}
\ee
\end{definition}

The existence and uniqueness of the weak solution of \refe{stocont} in the sense of the above definition is a classical result (again, see e.g. \cite{temam} or \cite{bf}).

\subsection{The finite volume scheme}

Under hypotheses \refe{hypomega}-\refe{hypfg}, let $\disc$ be an
admissible discretization of $\O$ in the sense of Definition \ref{adisc}.
It is then natural to write an approximate problem to the Stokes problem \refe{stocontf} in the following way.
\be
\left\{ \ba \dsp
u \in E_\disc(\O), \ p\in H_\disc(\O) \hbox{ with } \int_\O p(x)\d x = 0
\\[2ex] \dsp
\eta \int_\O  u(x)\cdot  v(x) \d x +  \nu [u, v]_{\disc}
\\ \dsp \hspace{7ex}
- \int_\O p(x) \div_\disc(v)(x) \d x =  \int_\O f(x)\cdot v(x) \d x
& \dsp \quad \forall v \in H_\disc(\O)^d
\ea\right.\label{schvf0}
\ee
As we use a colocated approximation for the velocity and the pressure fields, the scheme must be stabilized.
Using a non-consistant stabilization {\it \`a la} Brezzi-Pitk\"aranta \cite{brez}, we  then look for $(u,p)$ such that
\be
\left\{ \ba \dsp
(u,p)\in H_\disc(\O)^d\times H_\disc(\O) \hbox{ with } \int_\O p(x) \d x = 0
\\[2ex] \dsp
\eta \int_\O u(x)\cdot  v(x) \d x + \nu [u, v]_{\disc}
\\ \dsp \hspace{7ex}
- \int_\O p(x) \div_{\disc}(v)(x) \d x =  \int_\O  f(x)\cdot  v(x) \d x \qquad
& \dsp
\forall v \in  H_\disc(\O)^d
\\[3ex] \dsp
\int_\O \div_\disc(u)(x) q(x) \d x  = - \lambda \ \size(\disc)^\alpha\ \langle p,q \rangle_{\disc}
& \dsp
\forall q\in H_\disc(\O)
\ea \right.
\label{schvf}
\ee
where $\lambda >0$ and $\alpha\in(0,2)$ are adjustable parameters of the scheme which will have to be tuned in order to make a balance between accuracy and stability.

\medskip
System \refe{schvf} is equivalent to finding the family of vectors $(u_K)_{K\in\mesh} \subset \R^d$, and scalars $(p_K)_{K\in\mesh}\subset \R$ solution of the system of equations obtained by writing for each control volume $K$ of $\mesh$:
\be
\left\{ \ba \dsp
\eta\ \mcv\ u_\cv
- \nu \sum_{L\in\NN_K} \tkl (u_L - u_K)
- \nu \sum_{\s\in\edgescv\cap\edgesext} \tks (0 - u_K) \qquad
\\ \dsp \hfill
+ \sum_{L\in\NN_K} A_{KL}\ (p_L - p_K)  = \int_K f(x) \d x
\\[4ex]\dsp
\sum_{L\in\NN_K}   A_{KL}\cdot (u_K+u_L) - \lambda
\ \size(\disc)^\alpha \sum_{L\in\NN_K} \tkl (p_L - p_K) = 0
\ea \right.
\label{schvfS}
\ee
supplemented by the relation
\be
\sum_{K\in\mesh} \mK\ p_K = 0
\label{moypnulle}
\ee
Defining  $p_\edge = (p_K+p_L)/2 $ if $\edge = \edgecvcvv$, and $p_\edge = 
p_K$ if $\edge \in \edgesext \cap \edgescv$, and using the fact 
that $\sum_{\edge\in \edgescv} \medge \n_{K,\sigma} =0$, one notices that:
$\sum_{L\in\NN_K} A_{KL}\ (p_L - p_K)$ is in fact equal to 
$\sum_{\edge\in \edgescv} \medge p_\edge \n_{K,\sigma}, $ thus yielding a conservative form, 
which shows that \refe{schvfS} is indeed a finite volume scheme. 

\medskip
The existence of a solution to \refe{schvf} will be proven below.

\subsection{Study of the scheme in the linear case}\label{seccvglin}

We first prove a stability estimate for the velocity.

\begin{proposition}[Discrete $H^1$ estimate on velocities]\label{estl1l2}
Under hypotheses \refe{hypomega}-\refe{hypfg}, let $\disc$ be an admissible discretization of $\O$ in the sense of definition \ref{adisc}.
Let $\lambda\in(0,+\infty)$ and $\alpha\in(0,2)$ be given.
Let $(u,p)\in H_\disc(\O)^d\times H_\disc(\O)$ be a solution to \refe{schvf}.
Then the following inequalities hold:
\be
 \nu \Vert u\Vert_\disc \le {\rm diam}(\O) \Vert f \Vert_{(L^2(\O))^d},
\label{estimu}\ee
and
\be
 \nu\ \lambda \ \size(\disc)^\alpha \ \vert p \vert_\disc^2
\le{\rm diam}(\O)^2 \Vert f \Vert_{(L^2(\O))^d}^2.
\label{estimp}
\ee
\end{proposition}

\begin{proof} We apply \refe{schvf} setting $v = u$.
We get
\[
 \eta\int_\O u(x)^2 \d x+ 
 \nu \Vert u\Vert_\disc^2 - \int_\O p(x) \div_\disc(u)(x) \d x
= \int_\O  f(x)\cdot v(x) \d x.
\]
 Since $\eta \ge 0,$  the second equation of \refe{schvf} with $q=p$ and Young's inequality yield that:
\[
\ba\dsp
\eta \int_\O u(x)^2 \d x +
 \nu \Vert u\Vert_\disc^2 + \lambda \ \size(\disc)^\alpha\ \vert p\vert_\disc^2
  \le\\ \dsp
  \frac {{\rm diam}(\O)^2} {2\nu} \Vert f \Vert_{(L^2(\O))^d}^2
 + \frac {\nu} {2{\rm diam}(\O)^2} \Vert u  \Vert_{(L^2(\O))^d}^2.
\ea
\]
 Using the Poincar\'e inequality \refe{poindis} gives
\[
 \nu \Vert u\Vert_\disc^2 + \lambda \ \size(\disc)^\alpha\ \vert p\vert_\disc^2
 \le \frac {{\rm diam}(\O)^2} {2\nu} \Vert f \Vert_{(L^2(\O))^d}^2
 + \frac {\nu} {2} \Vert u\Vert_\disc^2,
\]
which leads to \refe{estimu} and \refe{estimp}.
%% RE ici il n'y avait pas a reporter l'une pour avoir l'autre
\end{proof}

\medskip
We can now state the existence and the uniqueness of a discrete solution to \refe{schvf}.
\begin{corollary}\label{exunpos}{\bf [Existence and uniqueness of a solution to the finite volume scheme]}
Under hypotheses \refe{hypomega}-\refe{hypfg}, let $\disc$ be an admissible discretization of $\O$ in the sense of Definition \ref{adisc}.
Let $\lambda\in(0,+\infty)$  and $\alpha\in(0,2)$ be given.
Then there exists a unique solution to \refe{schvf}.
\end{corollary}
\begin{proof}
System \refe{schvf} is a linear system.
Assume that $f=0$.
From propositions \ref{estl1l2} and using \refe{poindismoy}, we get that $u = 0$ and $p = 0$. This proves that the linear system \refe{schvf} is invertible.
\end{proof}

We then prove the following strong estimate on the pressures.

\begin{proposition}[$L^2$ estimate on pressures]\label{estp}
Under hypotheses \refe{hypomega}-\refe{hypfg}, let $\disc$ be an 
admissible discretization of $\O$ in the sense of definition \ref{adisc} and let $\theta >0$ be such that $\regul(\disc)> \theta$.
Let $\lambda\in(0,+\infty)$ and $\alpha\in(0,2)$ be given.
Let $(u,p)\in H_\disc(\O)^d\times H_\disc(\O)$ be a solution to \refe{schvf}.
Then there exists $\ctel{estisp}$, only depending on $d$, $\O$, $\eta$, $\nu$, $\lambda$, $\alpha$  and $\theta$, and not on $\size(\disc)$, such that the following inequality holds:
\be
\Vert p\Vert_{L^2(\O)} \le \cter{estisp} \Vert f \Vert_{(L^2(\O))^d}.
\label{estimsp}
\ee
\end{proposition}

\begin{proof} 
We first apply a result by Ne\v cas \cite{necas}: thanks to $\int_\O p(x) \d x = 0$, there exists $\ctel{derham}>0$, which only depends on $d$ and $\O$, and $\bar v \in H^1_0(\O)^d$ such that $\div \bar v(x) = p(x)$ for a.e. $x\in\O$ and
\be
\Vert \bar v\Vert_{H^1_0(\O)^d}\le \cter{derham}\Vert p\Vert_{L^2(\O)}.
\label{ineqderham}
\ee
We then set
\[
v_\s^{(i)} = \frac 1 {\medge} \int_\s \bar v^{(i)}(x) \d\gamma(x), 
\ \forall \sigma\in\edges,\ \forall i=1,\ldots,d.
\]
(note that $v_\sigma^{(i)} = 0$ for all $\sigma\in\edgesext$ and $i=1,\ldots,d$) and we define $v\in H_{\disc}(\O)^d$ by
\[
v_K^{(i)} = \frac 1 {\mK} \int_K \bar v^{(i)}(x) \d x, \ \forall 
K\in\mesh,\ \forall i=1,\ldots,d.
\]
Applying the results given p 777 in \cite{book}, we get that there exists $\ctel{dirnh}>0$, only depending on $d$ and $\theta$, such that
\be
(v_K^{(i)} - v_\s^{(i)})^2 \le \cter{dirnh} \frac { \diam(K)} {\medge} 
\int_K (\grad v^{(i)}(x))^2 \d x,
\label{inebook1}
\ee
and 
\be
\Vert v\Vert_\disc\le \cter{dirnh} \Vert \bar v\Vert_{H^1_0(\O)^d}.
\label{continjh1d}
\ee
We then have
\[
\int_\O p(x) \div_\disc v(x) \d x = \sum_{K\in\mesh} p_K 
\sum_{L\in\NN_K}   A_{KL} \cdot (v_K+v_L) = \terml{esp1}+ \terml{esp2},
\]
where
\[
\ba\dsp
\termr{esp1} &=&\dsp \sum_{K\in\mesh} p_K 
\sum_{L\in\NN_K}    2 A _{KL} \cdot v_{K|L} \\ &=& \dsp 
\sum_{K\in\mesh} p_K \sum_{L\in\NN_K} \int_{K|L} \bar v(x)\cdot \n_{KL}  
\d\gamma(x) \\ &=& \dsp
\int_\O p(x) \div \bar v(x) \d x = \Vert p\Vert_{L^2(\O)}^2,
\ea
\]
and
\[
\ba\dsp
\termr{esp2} &=&\dsp \sum_{K\in\mesh} p_K 
\sum_{L\in\NN_K}  \mkl\left(\frac 1 2(v_K+v_L) -  v_{K|L}\right)\cdot 
\n_{KL}\\ &=& \dsp 
\sum_{\s=K|L\in\edgesint} \mkl (p_K - p_L) \left(\frac 1 2(v_K+v_L) -  
v_{K|L}\right)\cdot \n_{KL}.
\ea
\]
We then have, thanks to the Cauchy-Schwarz inequality
\[
\termr{esp2}^2 \le \vert p\vert_\disc^2  \sum_{\s=K|L\in\edgesint} 
 \mkl d_{KL}\left(\frac 1 2(v_K+v_L) -  v_{K|L}\right)^2.
\]
Applying Inequality \refe{inebook1} and thanks to
$(\frac 1 2(v_K+v_L)-v_{K|L})^2 \le \frac 1 2 ((v_K -  v_{K|L})^2 + (v_L -  
v_{K|L})^2)$, we get that
\[
\termr{esp2}^2 \le \vert p\vert_\disc^2  \sum_{\s=K|L\in\edgesint} 
d_{KL} \cter{dirnh} \size(\disc)  
\int_{K\cup L} \sumi (\grad v^{(i)}(x))^2 \d x.
\]
This in turn implies the existence of $\ctel{estisp1}>0$, only depending on $d$ and $\theta$, such that
\[
\termr{esp2}^2 \le 
\cter{estisp1}\size(\disc)^2 \vert p\vert_\disc^2  \Vert \bar 
v\Vert_{H^1_0(\O)^d}^2.
\]
Thanks to \refe{ineqderham}, we then get, gathering the previous results
\be
\int_\O p(x) \div_\disc v(x) \d x \ge \Vert p\Vert_{L^2(\O)}^2 - 
\cter{estisp1}\size(\disc) \vert p\vert_\disc 
\cter{derham}\Vert p\Vert_{L^2(\O)}.
\label{latche}
\ee
We then introduce $v$ as a test function in \refe{schvf}.
We get
\be
\int_\O p(x) \div_{\disc}(v)(x) \d x = \eta \int_\O u(x)\cdot  v(x) \d x +
 \dsp \nu [u, v]_{\disc} - \int_\O  f(x)\cdot  v(x)\d x.
\label{vdansschvf}
\ee
Applying the discrete Poincar\'e inequality, \refe{continjh1d} and 
\refe{latche}, we get the existence of  $\ctel{estisp3}$, only depending on
$d$, $\O$, $f$, $\eta$, $\nu$, $\lambda$ and $\theta$, such that
\[
\Vert p\Vert_{L^2(\O)}^2 - 
\cter{estisp1}\size(\disc) \vert p\vert_\disc
\cter{derham}\Vert p\Vert_{L^2(\O)}
\le \cter{estisp3} \left( \Vert u\Vert_\disc + \Vert f\Vert_{L^2(\O)^d}\right) 
\Vert p\Vert_{L^2(\O)}.
\]
We now apply \refe{estimu} and \refe{estimp}.
Since $\size(\disc)^2 \le \size(\disc)^\alpha \diam(\O)^{2-\alpha}$, the condition
$\alpha\le 2$ suffices to produce \refe{estimsp} from the above inequality, a factor
$1/\lambda$ being introduced in the expression of  $\cter{estisp}$ (it is therefore not possible to let $\lambda$ tend to 0  in \refe{estimsp}).
\end{proof}

We then have the following result, which states the convergence of the
scheme \refe{schvf}.

\begin{proposition}[Convergence  in the linear case]\label{cvgce}
Under hypotheses \refe{hypomega}-\refe{hypfg}, let $(\bar u,\bar p)$ be the unique  weak solution of the Stokes problem \refe{stocont} in the sense of definition \ref{weaksol}.
Let $\lambda\in(0,+\infty)$, $\alpha\in(0,2)$ and $\theta>0$ be given and let $\disc$
be an admissible discretization of $\O$ in the sense of definition \ref{adisc} 
such that $\regul(\disc)\ge \theta$.
Let  $(u,p)\in H_{\disc}(\O)^d \times H_{\disc}(\O)$ be the unique solution to \refe{schvf}.

Then $u$ converges to $\bar u$ in $(L^2(\O))^d$ and $p$ weakly converges to $\bar p$ in $L^2(\O)$ as $\size(\disc)$ tends to $0$.
\end{proposition}

\begin{proof}
Under the hypotheses of the above proposition, let $(\disc^{(m)})_{m\in\N}$ be a sequence of admissible discretizations of $\O$ in the sense of definition \ref{adisc}, such that $\lim_{m\to\infty}\size(\disc^{(m)}) = 0$ and such that $\regul(\disc^{(m)})\ge \theta$, for all $m\in\N$.
\\
Let $(u^{(m)},p^{(m)})\in H_{\disc^{(m)}}(\O)^d \times H_{\disc^{(m)}}(\O)$ be given
by \refe{schvf} for all $m\in\N$.
Let us prove the existence of a subsequence of $(\disc^{(m)})_{m\in\N}$ such 
that the corresponding sequence  $(u^{(m)})_{m\in\N}$ converges in $(L^2(\O))^2$ 
to $\bar u$ and the sequence $(p^{(m)})_{m\in\N}$ weakly converges in $(L^2(\O))^2$ to $\bar p$, as $m\to\infty$.
Then the proof is complete thanks to the uniqueness of $(\bar u,\bar p)$.

Using \refe{estimu}, we obtain (see \cite{cvnl}, \cite{book}) an estimate on the translates of $u^{(m)}$:  for all $m\in\N$, there exists $\ctel{cc2}>0$, only depending on $\O$, $\nu$, $f$ and $g$ such that
\be
\ba
\int_\O(u^{(m,k)}(x+\xi)-u^{(m,k)}(x))^2 \d x \le \cter{cc2}
|\xi|(|\xi| + 4\size(\disc^{(m)})), \\
\mbox{ for } k = 1,\ldots,d,\ \forall 
\xi\in\R^d,
\ea \label{transx}
\ee
where $u^{(m,k)}$ denotes the $k$-th component of $u^{(m)}$.
We may then  apply Kolmogorov's theorem, and obtain the existence of a subsequence of  $(\disc^{(m)})_{m\in\N}$ and of $\bar u \in H^1_0(\O)^2$ such that $(u^{(m)})_{m\in\N}$
converges to $\bar u$ in $L^2(\O)^2$.
Thanks to proposition \ref{estp}, we extract from this subsequence another 
one (still denoted $u^{(m)}$) such that $(p^{(m)})_{m\in\N}$ weakly 
converges to some function $\bar p$ in $L^2(\O)$.
In order to  conclude the proof of the convergence of the scheme, 
there only remains to   prove that $(\bar u,\bar p)$ 
is the solution of \refe{stocontf}, thanks to the uniqueness of this solution.

Let $\varphi\in (C^\infty_c(\O))^d$.
Let $m\in\N$ such that $\disc^{(m)}$ belongs to the above extracted subsequence and let
$(u^{(m)},p^{(m)})$  be the solution to \refe{schvf} with $\disc = \disc^{(m)}$.
We suppose that $m$ is large enough and thus $\size(\disc^{(m)})$ is small 
enough to ensure for all $K\in\mesh$ such that
 $K\cap$ support$(\varphi)\ \neq\emptyset$, then $\dr K\cap \dr\O = \emptyset$ 
 holds.
Let us take $v = P_{\disc^{(m)}}\varphi$ in \refe{schvf}.
Applying proposition \ref{cpct}, we get
\[
\lim_{n\to\infty} [u^{(m)},P_{\disc^{(m)}}\varphi]_{\disc^{(m)}} =
\int_\O \nabla \bar u(x):\nabla\varphi(x) \d x.
\]
Moreover, it is clear that
\[
\lim_{n\to\infty} \int_\O f(x)\cdot P_{\disc^{(m)}}\varphi(x)\d x =
 \int_\O f(x)\cdot \varphi(x)\d x,
\]
and 
\[
\lim_{n\to\infty}\eta \int_\O u^{(m)}(x) \cdot  P_{\disc^{(m)}}\varphi(x)\d x =
\eta \int_\O \bar u(x)\cdot\varphi(x)\d x.
\]
Thanks to the weak convergence of the sequence of approximate pressures, to \refe{estimp}
and to the hypothesis $\alpha < 2$, we now apply proposition \ref{weakconvgrad}, which gives
\be
\lim_{n\to\infty}
\int_\O p^{(m)}(x)\div_{\disc^{(m)}}(P_{\disc^{(m)}}\varphi)(x)
\d x = \int_\O \bar p(x)\div\varphi(x)
\d x.
\label{convp}
\ee

\medskip
The last step is to prove that $\div(\bar u) = 0$ a.e. in $\O$.
Let $\varphi\in C^\infty_c(\O)$ and let $m\in\N$ be given.
Let us take $q = P_{\disc^{(m)}}\varphi$ in \refe{schvf}.
We get $\terml{X}^{(m)} = -\terml{Y}^{(m)}$, where
\[
\termr{X}^{(m)} = \int_\O \div_{\disc^{(m)}}(x) (u^{(m)}) P_{\disc^{(m)}}\varphi(x) \d x.
\]
and
\[
\termr{Y}^{(m)} = \lambda\ \size(\disc^{(m)})^{\alpha}
\langle p^{(m)},P_{\disc^{(m)}}\varphi\rangle_\disc.
\]
On the one hand, the third item of proposition \ref{cpct} produces
\[
\lim_{n\to\infty}\termr{X}^{(m)} =  \sumi \int_\O \varphi(x) \partial_i 
\bar u^{(i)} \d x.
\]
On the other hand, using the Cauchy-Schwarz inequality, we get:
\[
 \termr{Y}^{(m)}  \le  \lambda \size(\disc^{(m)})^{\alpha}  
\vert  p^{(m)}\vert_\disc 
\vert P_{\disc^{(m)}}\varphi\vert_\disc 
\]
Therefore, thanks to \refe{estimp} and to the regularity of $\varphi$ (that implies that
$\vert P_{\disc^{(m)}}\varphi\vert_\disc$ remains bounded independently on $\size(\disc^{(m)})$) we obtain $\lim_{n\to\infty} \termr{Y}^{(m)} = 0$.
This in turn implies that:
\be
\sumi \int_\O \varphi(x) \partial_i \bar u^{(i)}(x) \d x = 0,
\mbox{ for all  }\varphi\in C^\infty_c(\O), 
\label{vandiv}
\ee
which proves that $\bar u \in E(\O)$.
\end{proof}

\begin{remark}[Strong convergence of the pressure] Note that the proof of 
the strong convergence of  $p$ to $\bar p$ is a straightforward consequence
of the error estimate stated in Proposition \ref{ester} below, which holds
under additional regularity hypotheses.
\end{remark}

\subsection{An error estimate}

We then have the following result, which states an error estimate for the 
scheme \refe{schvf}.

\begin{proposition}[Error estimate   in the linear case]\label{ester}
 Under hypotheses \refe{hypomega}-\refe{hypfg}, we assume that the weak solution  $(\bar u,\bar p)$ of the Stokes problem \refe{stocont} in the sense of definition \refe{weaksol} is such that $(\bar u,\bar p)\in H^2(\O)^d \times H^1(\O)$.
Let $\lambda\in(0,+\infty)$  and $\alpha\in(0,2)$ be given, let $\disc$ be an admissible
discretization of $\O$ in the sense of definition \ref{adisc} and let $\theta>0$ such that $\regul(\disc^{(m)})\ge \theta$.
Let $(u,p)\in H_{\disc}(\O)^d \times H_{\disc}(\O)$ be the solution to \refe{schvf}.
Then there exists $\ctel{ester1}$, which only depends on $d$, $\O$, $\nu$, $\eta$ and $\theta$ such that
\be
\Vert u - \bar u\Vert_{L^2(\O)}^2\le \cter{ester1}\eps(\lambda,\size(\disc),\bar p, \bar u),
\label{eqester1}
\ee
\be
\lambda \ \size(\disc)^\alpha\ \vert p\vert_\disc^2 \le \cter{ester1}
\eps(\lambda,\size(\disc),\bar p, \bar u) \label{eqester2}
\ee
\be
\Vert p - \bar p\Vert_{L^2(\O)}^2\le \cter{ester1}\eps(\lambda,\size(\disc),\bar p, \bar u).
\label{eqester3}
\ee
where 
\be\ba
\eps(\lambda,\size(\disc),\bar p, \bar u) = & \min\left(\lambda\size(\disc)^{\alpha},
\frac 1 {\lambda} \size(\disc)^{2-\alpha}\right) \\ &\times \left(
\Vert \bar p \Vert_{H^1(\O)}^2 +  \Vert \bar u \Vert_{H^2(\O)}^2\right).
\ea\label{epsilon}
\ee
\end{proposition}

\begin{proof}
We define  $(\hat u,\hat p)\in H_{\disc}(\O)^d \times H_{\disc}(\O)$ by  
$\hat u = P_{\disc} \bar u$, which means $\hat u_K = \bar u(x_k)$ 
for all $K\in\mesh$, and  $\hat p_K = \frac 1 {\mK} \int_K \bar p(x) \d x$ 
for all $K\in\mesh$.
Integrating the first equation of \refe{stocont} on $K\in\mesh$ gives
\be
\eta \int_K \bar u(x) \d x +\sum_{\s\in\edgescv} \left(
\ba-\nu \int_{\s} \grad 
\bar u(x) : \n_{K,\s} \d\gamma(x)+ \\
\int_{\s} \bar p(x) \n_{K,\s} \d\gamma(x)\ea\right) = \dsp\int_K f(x) \d x.
\label{stocontK}
\ee
We introduce, for $K \in \mesh$, 
$\eps^u_K = \hat u_K - \frac 1 {\mK} \int_K \bar u(x) \d x$,
and,  for $L\in\NN_K$:  \\
$R_{K,L} = \frac 1 {\dkl}  (\hat u_L - \hat u_K) - 
\frac 1 {\mkl} \int_{\s} \grad \bar u(x) : \n_{K,\s} \d\gamma(x)$, \\
and for $\s\in \edgescv\cap\edgesext$,
$R_{K,\s} = \frac 1 {\dcvedge} (0 - \hat u_K) -
\frac 1 {\medge}\int_{\s} \grad \bar u(x) : \n_{K,\s} \d\gamma(x)$; \\
moreover, we define for $L\in\NN_K$: 
$\eps^p_{K|L} = \frac 1 2 (\hat p_K+\hat p_L) - \frac 1 {\mkl}\int_{K|L}
\bar p(x) \d\gamma(x),$ 
and for $\s\in \edgescv\cap\edgesext$, 
$\eps^p_{\s} = \hat p_K - \frac 1 {\medge}\int_{\s} \bar p(x) \d\gamma(x)$.
Using these notations and the relation 
$\sum_{\s\in\edgescv} \medge \n_{K,\s} = 0$,
we get from \refe{stocontK}
\[
\ba
\eta\ \mcv \hat u_\cv  - \dsp
\nu\left(\sum_{L\in\NN_K} \tkl (\hat u_L - \hat u_K)  +  
\sum_{\s\in\edgescv\cap\edgesext} \tks (0 - \hat u_K)\right)+ \\ 
\dsp
\dsp \sum_{L\in\NN_K} A_{KL}\ (\hat p_L - \hat p_K)  =
 \dsp\int_K f(x) \d x + R_K,
\ea
\]
with
\[
\ba
R_K = \eta\ \mcv \eps^u_K   
- \dsp
\nu\left(\sum_{L\in\NN_K} {\mkl}R_{K,L}  +  
\sum_{\s\in\edgescv\cap\edgesext} \medge R_{K,\s} \right)+ 
\dsp \sum_{\s\in\edgescv} \medge \ \eps^p_{\s}\n_{K,\s}.
\ea
\]
We then set $\delta\!u = \hat u - u$ and $\delta\!p = \hat p - p$.
We then get, substracting the first relation of the scheme \refe{schvfS} to the above equation,
\be
\ba\eta \int_\O \delta\!u(x) v(x) \d x +
 \dsp \nu [\delta\!u,v]_{\disc}
- \int_\O \delta\!p(x) \div_{\disc}(v)(x) \d x = \\
 \int_\O R(x) v \d x,\ \forall v \in  H_\disc(\O)^d,
\ea \label{eqester4}
\ee
and, setting $v = \delta\!u$ in \refe{eqester4},
\[
\eta \int_\O \delta\!u(x)^2 \d x +
 \dsp \nu \Vert \delta\!u\Vert_{\disc}^2
- \int_\O \delta\!p(x) \div_{\disc}(\delta\!u)(x) \d x = 
 \int_\O R(x) \delta\!u(x) \d x.
\]
We now integrate the second equation of \refe{stocont} on $K\in\mesh$.
This gives
\[
\sum_{\s\in\edgescv} \int_{\s} \bar u(x)\cdot \n_{K,\s} \d\gamma(x) = 0,\ \forall K\in\mesh.
\]
Using $\bar u\in H^1_0(\O)$, we then obtain
\[
\sum_{L\in\NN_K}   A_{KL}\cdot (\hat u_K+\hat u_L) = 
\sum_{L\in\NN_K}   \mkl \eps^u_{K|L},\ \forall K\in\mesh
\]
with 
\[
\eps^u_{K|L} = \left(\frac 1 2 (\hat u_K+\hat u_L) - \frac 1 {\mkl} 
\int_{K|L} \bar u(x) \d\gamma(x) \right)\cdot \n_{KL},\ \forall K\in\mesh,\ \forall L\in\NN_K.
\]
We then give, substracting the second relation of the scheme  \refe{schvfS} to the above equation,
\[
\int_\O \div_\disc(\delta\!u)(x) \delta\!p(x) \d x  = 
\lambda \ \size(\disc)^\alpha\   \langle p,\hat p - p\rangle_{\disc}
+ \terml{esterm1}, \ 
\]
with
\[
\termr{esterm1} = \sum_{K|L\in \edgesint}   \mkl \eps^u_{K|L} (\delta\!p_K - 
\delta\!p_L),
\]
Gathering the above results, we get
\be
\ba\eta \int_\O \delta\!u(x)^2 \d x +
 \dsp \nu \Vert \delta\!u\Vert_{\disc}^2 + \lambda \ 
\size(\disc)^\alpha\   \vert p\vert_\disc^2 = \\
\lambda \ \size(\disc)^\alpha\   \langle p,\hat p\rangle_{\disc} + 
 \int_\O R(x)\cdot  \delta\!u(x) \d x + \termr{esterm1}.
\label{eqester5}
\ea\ee
Let us study the terms at the right hand side of the above equation.
We have, using the Young inequality,
\be
\langle p,\hat p\rangle_{\disc} \le \frac 1 4 \vert p\vert_\disc^2
+
\vert \hat p\vert_{\disc}^2 \le 
\frac 1 4 \vert p\vert_\disc^2 + \ctel{pbar} \Vert \bar p 
\Vert_{H^1(\O)}^2.
\label{eqester6}
\ee
We then study
$ \int_\O R(x)\cdot \delta\!u(x) \d x = \terml{esterm2}+\terml{esterm3}+\terml{esterm4}$,
with
\[
\termr{esterm2} =  \eta\ \int_\O \eps^u(x)\cdot \delta\!u(x) \d x,
\]
\[
\termr{esterm3} = \nu\sum_{K\in\mesh} \left(\sum_{L\in\NN_K} {\mkl} R_{K,L}  +  
\sum_{\s\in\edgescv\cap\edgesext} \medge R_{K,\s} \right) \cdot 
\delta\!u_K,
\]
and
\[
\termr{esterm4} = \dsp \sum_{K\in\mesh} \sum_{\s\in\edgescv} 
\medge \ \eps^p_{\s}\n_{K,\s} \cdot\delta\!u_K.
\]
Thanks to interpolation results proven in \cite{book} and to \refe{poindis}, we obtain
\be
\termr{esterm2}\le 
\ctel{ester2b} \size(\disc)^2 \Vert \bar u \Vert_{H^2(\O)}^2+ 
\frac \nu 4 \Vert \delta\!u \Vert_\disc^2,
\label{eqester7}
\ee
\be
\termr{esterm3} \le \ctel{ester3} \size(\disc)^2 \Vert \bar u 
\Vert_{H^2(\O)}^2+ 
\frac \nu 4 \Vert \delta\!u\Vert_\disc^2,
\label{eqester8}
\ee
and
\be
\termr{esterm4} \le \ctel{ester4} \size(\disc)^2 \Vert \bar p 
\Vert_{H^1(\O)}^2+ 
\frac \nu 4 \Vert \delta\!u\Vert_\disc^2.
\label{eqester9}
\ee
We then study $\termr{esterm1}$. We have $\termr{esterm1} = \terml{esterm1a}-\terml{esterm1b}$ with
\[
\termr{esterm1a} = \sum_{K|L\in \edgesint}   \mkl \eps^u_{K|L} (\hat p_K - 
\hat p_L),
\]
which verifies
\be
\termr{esterm1a} \le \ctel{ester1a} \size(\disc) \left(
\Vert \bar p \Vert_{H^1(\O)}^2 +  \Vert \bar u \Vert_{H^2(\O)}^2\right),
\label{eqester10}
\ee
and
\[
\termr{esterm1b} = \sum_{K|L\in \edgesint}   \mkl \eps^u_{K|L} (p_K - p_L),
\]
which verifies
\be
\termr{esterm1b} \le \frac 1 4 \lambda \ 
\size(\disc)^\alpha\   \vert p\vert_\disc^2
+\ctel{ester1b} \frac 1 {\lambda}  \size(\disc)^{2-\alpha}\Vert \bar u 
\Vert_{H^2(\O)}^2 .
\label{eqester10bis}
\ee
Gathering equations \refe{eqester5}-\refe{eqester10bis} gives
\[
\Vert \delta\!u\Vert_\disc^2 + \lambda \ 
\size(\disc)^\alpha\   \vert p\vert_\disc^2
\le \ctel{estertot}\eps(\lambda,\size(\disc),\bar p, \bar u),
\]
where $\eps(\lambda,\size(\disc),\bar p, \bar u)$ is defined by \refe{epsilon}. 
This in turn yields \refe{eqester1} and \refe{eqester2}.
We then again follow the method used in the proof of Proposition \ref{estp}.
Using $\int_\O \hat p(x) \d x = 0$ and therefore $\int_\O \delta\!p(x) \d x = 0$, let $\bar v \in H^1_0(\O)^d$ be given such that $\div \bar v(x) = \delta\!p(x)$ for a.e. $x\in\O$ and
\be
\Vert \bar v\Vert_{H^1_0(\O)^d}\le \cter{derham}\Vert
\delta\!p\Vert_{L^2(\O)}.
\label{eqester11}
\ee
We again set
\[
v_\s^{(i)} = \frac 1 {\medge} \int_\s \bar v^{(i)}(x) \d\gamma(x), 
\ \forall \sigma\in\edges,\ \forall i=1,\ldots,d.
\]
and we define $v\in H_{\disc}(\O)^d$ by
\[
v_K^{(i)} = \frac 1 {\mK} \int_K \bar v^{(i)}(x) \d x, \ \forall 
K\in\mesh,\ \forall i=1,\ldots,d.
\]
The same method gives
\[
\ba
\Vert \delta\!p\Vert_{L^2(\O)}^2 &\dsp \le \int_\O \delta\!p(x) 
\div_\disc(v)(x) \d x +
\cter{estisp1}\size(\disc)\vert p\vert_\disc   \Vert \bar 
v\Vert_{H^1_0(\O)^d} \\
&\dsp \le \int_\O \delta\!p(x) \div_\disc(v)(x) \d x +
\ctel{ester15}\size(\disc)^2\vert p\vert_\disc^2 + \frac 1 4 \Vert 
\delta\!p\Vert_{L^2(\O)}^2.
\ea
\]
We now use $v$ as test function in \refe{eqester4}.
We get
\[
\int_\O \delta\!p(x) \div_{\disc}(v)(x) \d x = \eta \int_\O \delta\!u(x) 
v(x) \d x +
 \dsp \nu [\delta\!u,v]_{\disc} + \int_\O R(x) v \d x.
\]
Gathering the two above inequalities, \refe{eqester7}, \refe{eqester8}, \refe{eqester9} and \refe{eqester11} produces
\[
\ba
\Vert \delta\!p\Vert_{L^2(\O)}^2 \le & \dsp  \frac 1 2 \Vert 
\delta\!p\Vert_{L^2(\O)}^2 + 
\ctel{ester17} \size(\disc)^2\left(
\Vert \bar p \Vert_{H^1(\O)}^2 +  \Vert \bar u \Vert_{H^2(\O)}^2\right) \\
&\dsp + \ctel{ester16} \Vert \delta\!u\Vert_\disc^2 + 
\cter{ester15}\size(\disc)^2\vert p\vert_\disc^2.
\ea
\]
Applying \refe{eqester1} and \refe{eqester2} gives \refe{eqester3}.
\end{proof}

%% RE une remarque sur .5, car on ne l'avait meme pas indique...
\begin{remark} In the above result, it suffices to let $\alpha = 1$ to obtain the proof of an
order $1/2$ for the convergence of the scheme. We recall that this result is not sharp, and
that the numerical results show a much better order of convergence.
\end{remark}

%-------------------------------------------------------------------------------------------
% Navier-Stokes
%-------------------------------------------------------------------------------------------

\section{The finite volume scheme for the Navier-Stokes equations}\label{secfvsnlin}

Before handling the transient nonlinear case, we first address in the following
 section the steady-state case.

\subsection{The steady-state case}\label{stnl}

For the following continuous equations,
\be
\ba\dsp
\eta  \bar u^{(i)}  -\nu \lap \bar u^{(i)} + \partial_i \bar p + \sumj \bar u^{(j)}
\partial_j \bar u^{(i)}
 = f^{(i)} 
\hbox{ in } \Omega, \ \mbox{ for } i = 1,\ldots,d,\\ \dsp
\div \bar u = \sumi \partial_i \bar u^{(i)}   = 0 \hbox{ in } \Omega.
\ea\label{nstocontss}
\ee
with a homogeneous Dirichlet boundary condition, we define the following weak sense.

\medskip
\begin{definition}[Weak solution for the steady Navier-Stokes equations]\label{weaksolss}
Under hypotheses \refe{hypomega}-\refe{hypfg}, let $E(\O)$ be defined by \refe{e0}.
Then  $(\bar u,\bar p) $ is called a weak solution of \refe{nstocontss} if
\be
\left\{\ba \dsp
\bar u \in E(\O), \ \bar p\in L^2(\O) \hbox{ with } \int_\O \bar p(x) \d x = 0,
\\[2ex]
\dsp \eta \int_\O \bar u(x)\cdot \bar v(x) \d x + \dsp\nu \int_\O \nabla \bar u(x) :\nabla \bar v(x) \d x
\\
\dsp \quad -\int_\O \bar p(x)\div \bar v(x) \d x  + b(\bar u,\bar u,\bar v)= \dsp \int_\O f(x)\cdot  \bar v(x) \d x
 \qquad  \forall \bar v \in H^1_0(\O)^d,
\ea\right.\label{nstocontfss}
\ee
where the trilinear form $b(.,.,.)$ is defined by \refe{deftricont}.
\end{definition}

\medskip
We now give the finite volume scheme for this problem.
Under hypotheses \refe{hypomega}-\refe{hypfg}, let $\disc$ be an admissible discretization of $\O$ in the sense of Definition \ref{adisc}.
We introduce Bernoulli's pressure $p + \half u^2$ instead of $p$, again denoted by $p$, and
for any real value $\lambda >0$ and $\alpha\in(0,2)$, we look for $(u,p)$ such that
\be
\qquad \left\{ \ba \dsp
(u,p)\in H_\disc(\O)^d\times H_\disc(\O) \hbox{ with } \int_\O p(x) \d x =0,
\\[4ex] \dsp
\eta \int_\O u(x)\cdot  v(x) \d x +
 \dsp \nu [u, v]_{\disc} + \half \int_\O u(x)^2 \div_{\disc}(v)(x)\d x
\\[2ex] \hspace{2ex} \dsp
- \int_\O p(x) \div_{\disc}(v)(x) \d x + b_\disc(u,u,v)=
 \int_\O  f(x)\cdot  v(x) \d x
& \dsp
 \forall v \in  H_\disc(\O)^d
\\[4ex] \dsp
\int_\O \div_\disc(u)(x) q(x) \d x  = - \lambda \ \size(\disc)^\alpha\
\langle p,q\rangle_{\disc}
& \dsp
\forall q\in H_\disc(\O)
\ea \right.
\label{schvfnlss}
\ee
where, for $u,v,w \in H_\disc(\O)$,  we define the following approximation for $b(u,v,w)$
\be
b_\disc(u,v,w) = \half \sum_{K\in\mesh} \sum_{L\in\NN_K}
(A_{KL} \cdot (u_K + u_L))\ ((v_L - v_K) \cdot w_K)
\label{deftridc}
\ee

\medskip
System \refe{schvfnlss} is equivalent to finding the family of vectors $(u_K)_{K\in\mesh} \subset \R^d$, and scalars $(p_K)_{K\in\mesh}\subset \R$ solution of the system of equations obtained by writing for each control volume $K$ of $\mesh$:
\be
\left\{ \ba \dsp
\eta\ \mcv\ u_\cv
- \nu \sum_{L\in\NN_K} \tkl (u_L - u_K)
- \nu \sum_{\s\in\edgescv\cap\edgesext} \tks (0 - u_K) \qquad
\\ \dsp
+ \sum_{L\in\NN_K} (A_{KL} \cdot (\half (u_K+u_L)))\ (u_L-u_K)
\\ \dsp 
+ \sum_{L\in\NN_K} A_{KL}\ (p_L - p_K)  
- \frac 1 2 \sum_{L\in\NN_K} A_{KL}\ (u_L^2 - u_K^2) = \int_K f(x) \d x
\\[4ex]\dsp
\sum_{L\in\NN_K}   A_{KL} \cdot (u_K+u_L) - \lambda
\ \size(\disc)^\alpha \sum_{L\in\NN_K} \tkl (p_L - p_K) = 0
\ea \right.
\label{schvfNSS}
\ee
supplemented by the relation:
\[
\sum_{K\in\mesh} \mK\ p_K = 0
\]

Defining $\tilde p_K = p_K  - u_K^2 /2$ and   
$\tilde p_\edge = (\tilde p_K + \tilde p_L)/2$ 
if $\edge = \edgecvcvv$,   $\tilde p_\edge = \tilde p_K$ 
if $\edge \in \edgesext \cap \edgescv$, and using the fact 
that $\sum_{\edge\in \edgescv} \medge \n_{K,\sigma} =0$, one again notices that:
$\sum_{L\in\NN_K} A_{KL}\ (\tilde p_L - \tilde p_K)$ is in fact equal to 
$\sum_{\edge\in \edgescv} \medge \tilde p_\edge \n_{K,\sigma}, $ 
thus yielding a conservative form for the fifth and sixth terms of the left 
handside of the discrete momentum equation in \refe{schvfNSS}. 
Defining $u_\edge = (u_K + u_L)/2$ 
if $\edge = \edgecvcvv$,   $u_\edge = 0$ 
if $\edge \in \edgesext \cap \edgescv$, one obtains that the nonlinear 
convective term  $\sum_{L\in\NN_K} (A_{KL} \cdot (\half (u_K+u_L)))\ (u_L-u_K)$ is equal to
$\sum_{\edge\in \edgescv} \medge   (\n_{K,\sigma} \cdot u_\edge) u_\edge
- \mcv u_\cv (\div_\disc u)_K $;  one may note that
 $(\div_\disc u)_K =  \sum_{\edge\in \edgescv} \medge   \n_{K,\sigma} \cdot u_\edge$.
Hence the nonlinear convective term is the sum of a conservative form and 
a source term due to the stabilization (this source term vanishes for a discrete divergence
free function $u$).

\medskip
Let us then study some properties of the trilinear form $b_\disc$. 
First note that the quantity $b_\disc(u,v,w)$ also writes
\be
b_\disc(u,v,w) = \half \sum_{K|L\in \edgesint}
(A_{KL} \cdot (u_K + u_L))\ ((v_L - v_K) \cdot (w_L + w_K))
\label{deftridcbis}
\ee
We thus get that, for all $u,v\in H_\disc(\O)^d$,
\be
\ba \dsp
b_\disc(u,v,v)
& \dsp
= \half \sum_{K|L\in \edgesint}
(A_{KL} \cdot  (u_K + u_L)) ((v_L)^2 - (v_K)^2)
\\[4ex] & \dsp
= -\half \int_\O  v(x)^2\ \div_{\disc}(u)(x)\ \d x
\ea
\label{superb}
\ee
We get in particular, that, for all $u\in E_\disc(\O)$, $b_\disc(u,u,u) = 0$, which is the discrete equivalent of the continuous property.

\medskip
\begin{remark}{\bf [Upstream weighting versions of the scheme]}
All the results of this paper are available, setting $F_{KL}(u) = A_{KL}\cdot(u_K + u_L)$ and considering, for $u,v,w\in H_\disc(\O)$,
\[
b_\disc^{\rm ups}(u,v,w) = b_\disc(u,v,w) + \frac 1 2
\sum_{K|L\in \edgesint} \Theta_{KL} | F_{KL}(u) |\  (v_L - v_K)\cdot (w_L - w_K),
\]
with, for example,  $ \Theta_{KL} = \max( 1 -  2 \nu \tkl/| F_{KL}(u) |, 0) $.
We then get, for all $u,v\in H_\disc(\O)$, the inequality
\[
b_\disc^{\rm ups}(u,v,v) \ge -\half \int_\O  v(x)^2\div_{\disc}(u)(x) \d x,
\]
which is sufficient to get all the estimates of this paper, 
together with the convergence properties of the scheme.
The use of such a local upwinding technique may be useful to avoid 
the development of  nonphysical oscillations only where meshes are too coarse.
\end{remark}

\medskip
The following technical estimates are crucial to prove the  
convergence properties of the scheme.

\medskip
\begin{lemma}[Estimates on $b_\disc(.,.,.)$ by discrete Sobolev norms]\label{l4l2h1} 
Under hypotheses \refe{hypomegat}-\refe{hypfgt}, let $\disc$ 
be an admissible discretization of $\O\times(0,T)$ in the 
sense of definition \ref{adisct}, and $\theta >0$ such that $\regul(\disc) \ge \theta$.
Then there exists $\ctel{bd4} >0$ and $\ctel{bd} >0$, only depending on 
$d$, $\theta$ and $\O$, such that
\be
b_\disc(u,v,w)\le\cter{bd4}\ \Vert u \Vert_{L^4(\O)^d}\ \Vert v\Vert_\disc\ \Vert w \Vert_{L^4(\O)^d}
\le \cter{bd}\ \Vert u\Vert_\disc\ \Vert v\Vert_\disc\ \Vert w\Vert_\disc.
\label{estib}
\ee
\end{lemma}

\begin{proof}
The quantity $b_\disc(u,v,w)$ reads
\[
\ba
b_\disc(u,v,w) =  \dsp\frac 1 4 \sum_{K\in\mesh}
\sum_{L\in\NN_K} (w_K \cdot  (v_L - v_K))
\ \tkl\ ((x_L - x_K)\cdot  (u_K + u_L))
\ea
\]
%% RE petite modif de forme, et yield au lieu de yields
Applying  the Cauchy-Schwarz inequality twice
and using the fact that $(x_L - x_K)^2 = d_{KL}^2$ and that, 
for any admissible discretization 
$\sum_{L\in\NN_K} \tkl\ d_{KL}^2 \le d\ \frac {m_K} \theta$ yield:
\[
\ba
b_\disc(u,v,w)^2
& \le & \dsp \ctel{dse}
& \dsp
\left(  \sum_{K\in\mesh} \sum_{L\in\NN_K} \tkl (w_K)^2 (x_L - x_K)^2 (2(u_K)^2 + 2(u_L)^2) \right)
\\[2ex]
& & & \dsp
\left( \sum_{K\in\mesh} \sum_{L\in\NN_K}  \tkl (v_L - v_K)^2 \right)
\\[2ex]
&\le &  \ctel{dse3}
& \dsp
\left( \sum_{K\in\mesh} m_K |w_K|^4  \right)^{1/2}
\left( \sum_{K\in\mesh}m_K |u_K|^4 \right)^{1/2} \Vert v\Vert_{\disc}^2.
\ea
\]
The inequality \refe{estib} is now a straightforward consequence 
of the following discrete Sobolev inequality, which holds under the same regularity
assumptions on the mesh (see proof in \cite{coudiere} or 
\cite[pp. 790-791]{book}):
\be
\dsp
\Vert u\Vert_{L^4(\O)} \le \ctel{l4hd}\ \Vert u\Vert_{\disc}.
\label{sob1}
\ee

\end{proof}

\begin{remark}[Two dimensional case]
In the case $d=2$, it may be proven setting $\alpha =2, p=p'=2$ in the proof 
p791 of \cite{book}, that \[
\Vert u\Vert_{L^4(\O)} \le \ctel{l2l4x}\Vert u\Vert_{L^2(\O)}^{1/2}
\Vert u\Vert_\disc^{1/2}
\]
and therefore, that 
there exists $\ctel{l2l4} >0$, only depending on $d$ and $\O$,  such that
\[
b_\disc(u,v,w) \le \cter{l2l4}\Vert v\Vert_\disc  \left(\Vert u\Vert_\disc\ \Vert
u\Vert_{L^2(\O)}\ \Vert w\Vert_\disc
\ \Vert w\Vert_{L^2(\O)} \right)^{1/2}.
\]
This is a discrete analogue to the classical continuous estimate on the
trilinear form. 
\end{remark}

%% RE sake au lieu de same
The existence of a solution to the scheme \refe{schvfnlss} is obtained 
through a so-called ``topological degree" argument. For the sake of 
completeness, we recall this argument  (which was first used for numerical
schemes in 
\cite{eggh}) 
in the finite dimensional case 
in the following theorem and refer to \cite{deimling} for the general case. 

\begin{theorem}[Application of the topological degree, finite dimensional case]
\label{degtop}
Let  $V$ be a finite dimensional vector space on $\R$ and 
$g$ be a continuous function from $V$ to $V$.  
Let us assume that there exists a continuous function $F$ 
from $V\times[0,1]$ to $V$ satisfying:
\begin{enumerate}
\item $F(\cdot,1) = g$, $F(\cdot,0) $ is an affine function.
\item There exists $R >0$, such that for any
$(v,\rho) \in V\times[0,1]$, if $F(v,\rho) = 0$, then  $\Vert v \Vert_V \ne R$.
\item The equation $F(v,0) = 0$ has a solution $v\in V$ such that
$\Vert v \Vert_V < R$.
\end{enumerate}

Then there exists at least a solution $v\in V$ such that $g(v) = 0$ and 
$\Vert v \Vert_V < R$.
\end{theorem}

Here $g(v) = 0$ represents the nonlinear system \refe{schvfnlss}, and 
we are now going to construct the function $F$ and show the required estimates.
Note that here, the use of Bernouilli's pressure leads to simpler calculations.

\medskip
\begin{proposition}[Discrete $H^1_0(\O)$ estimate on the velocities]\label{estl1l2ss}
Under hypotheses \refe{hypomega}-\refe{hypfg}, let $\disc$ 
be an admissible discretization of $\O\times(0,T)$ in the sense 
of definition \ref{adisct}.
Let $\lambda\in(0,+\infty)$  and $\alpha\in(0,2)$ be given.
Let $\rho\in[0,1]$ be given and let $(u,p)\in (H_\disc(\O))^d\times H_\disc(\O)$, 
be a solution to the following system of equations (which reduces to 
\refe{schvfnlss} as $\rho=1$ and to \refe{schvf} as $\rho = 0$)
%% RE y avait un pb de ref ... sch transitoire
\be
\qquad \left\{ \ba \dsp
(u,p)\in H_\disc(\O)^d\times H_\disc(\O) \hbox{ with } \int_\O p(x) \d x = 0,
\\[3ex] \dsp
\eta \int_\O u(x)\cdot  v(x) \d x +
 \dsp \nu [u, v]_{\disc} + \frac \rho 2  \int_\O u(x)^2 \div_{\disc}(v)(x)\d x
\\ \dsp
 + \rho\ b_\disc(u,u,v)
- \int_\O p(x) \div_{\disc}(v)(x) \d x =
 \int_\O  f(x)\cdot  v(x) \d x \hspace{1ex}
& \dsp
\forall v \in  H_\disc(\O)^d
\\[3ex] \dsp
\int_\O \div_\disc(u)(x) q(x) \d x  = - \lambda \ \size(\disc)^\alpha\
\langle p,q\rangle_{\disc}
& \dsp
\forall q\in H_\disc(\O)
\ea \right.
\label{schvfnlssrho}
\ee
Then $u$ and $p$ satisfy the following estimates, which are the same inequalities as obtained in the linear case (inequalities \refe{estimu} and  \refe{estimp}):
\[
\ba \dsp
 \nu \Vert u\Vert_\disc \le {\rm diam}(\O) \Vert f \Vert_{(L^2(\O))^d}
\\[2ex] \dsp
 \nu\ \lambda \ \size(\disc)^\alpha \ \vert p\vert_\disc^2
\le{\rm diam}(\O)^2  \Vert f \Vert_{(L^2(\O))^d}^2
\ea
\]
\end{proposition}

\begin{proof} 
The proof is similar to that of Proposition \ref{estl1l2}, 
using the property \refe{superb} on the discrete trilinear form.
\end{proof}

\medskip
We are now in position to prove the existence of at least one solution to scheme \refe{schvfnlss}.

\medskip
\begin{proposition}[Existence of a discrete solution]\label{exunss}
Under hypotheses \refe{hypomega}-\refe{hypfg}, let $\disc$ be an admissible discretization of $\O\times(0,T)$ in the sense of definition \ref{adisct}.
Let $\lambda\in(0,+\infty)$  and $\alpha\in(0,2)$ be given.
Then there exists at least one $(u,p)\in (H_\disc(\O))^d\times H_\disc(\O)$, solution to \refe{schvfnlss}.
\end{proposition}

\medskip
\begin{proof}
Let us define $V = \{(u,p)\in (H_\disc(\O))^d\times H_\disc(\O)$ 
 s.t. $\int_\O p(x) \d x = 0\}$. Consider the continuous application 
$F~:~V\times[0,1] 
\to V$ such that, for a given 
$(u,p)\in V$ 
  and $\rho\in[0,1]$, $(\hat u,\hat p) = F(u,p,\rho)$ 
 is defined by
\[
\ba \dsp
\int_\O \hat u(x)\cdot v(x) \d x
& \dsp
= \eta \int_\O u(x)\cdot  v(x) \d x +
\nu [u,v]_{\disc} -
\int_\O p(x) \div_{\disc}(v)(x) \d x
\\[3ex]  & \dsp
+ \rho \left(\half \int_\O u(x)^2 \div_{\disc}(v)(x) \d x
+  b_\disc(u,u,v) \right)
\\[3ex]  & \dsp
- \int_\O f(x) \cdot v(x)\d x
\hfill \forall  v \in H_\disc(\O)^d
\\[4ex] \dsp
\int_\O \hat p(x)\cdot q(x) \d x
& \dsp
= \int_\O \div_\disc(u)(x) q(x) \d x  +
\lambda \ \size(\disc)^\alpha\
\langle p,q\rangle_{\disc}
\hspace{7ex} \forall q\in H_\disc(\O).
\ea
\]
It is easily checked that the two above relations define a one to one function 
$F(.,.,.)$.
Indeed, the value of $\hat u^{(i)}_K$ and $\hat p_K$ for a given $K\in\mesh$ 
and $i=1,\ldots,d$ are readily obtained by setting $v^{(i)} = 1_K$, $v^{(j)} = 0$ 
for $j\neq i$,  and $q = 1_K$.

The application $F(.,.,.)$ is continuous, and, for a given $(u,p)$ such that
 $F(u,p,\rho) = (0,0)$, we can apply proposition \ref{estl1l2ss} and \refe{poindismoy}, 
 which prove that $(u,p)$ is bounded independently on $\rho$.
Since $F(u,p,0)$ is an affine function of $(u,p)$ 
(indeed invertible, see corollary \ref{exunpos}), 
we may apply Theorem \ref{degtop} 
and conclude to the existence of   at least one solution $(u,p)$ to \refe{schvfnlss}.
\end{proof}

\bigskip
We then have the following strong estimate on the pressures.

\medskip
\begin{proposition}[$L^2$ estimate on pressures]\label{estpss}
Under hypotheses \refe{hypomega}-\refe{hypfg}, let $\disc$ be an admissible
discretization of $\O$ in the sense of definition \ref{adisc}, and let $\theta>0$
such that $\regul(\disc)> \theta$.
Let $\lambda\in(0,+\infty)$ and $\alpha\in(0,2)$ be given.
Let $(u,p)\in H_\disc(\O)^d\times H_\disc(\O)$ be a solution to \refe{schvfnlss}.
Then there exists $\ctel{estispnlss}$, only depending on $d$, $\O$, $\eta$, $\nu$, $\lambda$, $\alpha$  and $\theta$, and not on $\size(\disc)$, such that the following inequality holds:
\be
\Vert p\Vert_{L^2(\O)} \le \cter{estispnlss} 
\left(\Vert f \Vert_{(L^2(\O))^d}+\left(\Vert f \Vert_{(L^2(\O))^d}\right)^2\right)
\label{estimspnlss}
\ee
\end{proposition}

\medskip
\begin{proof} 
We may follow  the proof of proposition \ref{estp}  until \refe{vdansschvf}, 
which is changed  to:
\be
\ba \dsp
\int_\O p(x) \div_{\disc}(v)(x) \d x
& \dsp
= \eta \int_\O u(x)\cdot  v(x) \d x +
\nu [u, v]_{\disc} - \int_\O  f(x)\cdot  v(x)\d x
\\ &\dsp
+ \half  \int_\O u(x)^2 \div_{\disc}(v)(x)\d x + b_\disc(u,u,v).
\ea
\label{vdansschvfnlss}
\ee
We again apply the discrete Poincar\'e inequality \refe{poindis}, \refe{continjh1d}, \refe{latche} and we use \refe{estib}.
We get the existence of  $\ctel{estispnl3}$, only depending on $d$, $\O$, $f$, $\eta$, $\nu$, $\lambda$ and $\theta$, such that
\[
\ba \dsp
\Vert p\Vert_{L^2(\O)}^2 - 
\cter{estisp1}\size(\disc) \vert p\vert_\disc
\cter{derham}\Vert p\Vert_{L^2(\O)} \le
\\[1ex] \dsp \hspace{20ex}
\cter{estispnl3} \left( \Vert u\Vert_\disc + \Vert f\Vert_{L^2(\O)^d} + \Vert u\Vert_\disc^2\right)
\Vert p\Vert_{L^2(\O)}
\ea
\]
We now apply \refe{estimu} and \refe{estimp}, which yields the conclusion.
\end{proof}

\bigskip
We now can state the convergence of Scheme \refe{schvfnlss}.

\medskip
\begin{theorem}[Convergence of the  scheme]\label{cvgcenlss}
Under hypotheses \refe{hypomega}-\refe{hypfg}, let $(\disc^{(m)})_{m\in\N}$ be a sequence of admissible discretizations of $\O$ in the sense of definition \ref{adisc}, such that
$\size(\disc^{(m)})$ tends to $0$ as $m\to\infty$ and such that there exists $\theta>0$ with $\regul(\disc^{(m)})\ge \theta$, for all $m\in\N$.
Let $\lambda\in(0,+\infty)$  and $\alpha\in(0,2)$ be given.
Let, for all $m\in\N$, $(u^{(m)},p^{(m)})\in (H_{\disc^{(m)}}(\O))^d\times H_{\disc^{(m)}}(\O)$, be a solution to 
%% RE \refe{schvfini}-\refe{solap} ici on avait laisse le schema du transitoire...
\refe{schvfnlss} with $\disc = \disc^{(m)}$.
Then there exists a weak solution $(\bar u,\bar p)$ of \refe{nstocontss} 
in the sense of definition \ref{weaksolss} and a subsequence of  $(\disc^{(m)})_{m\in\N}$, again denoted
$(\disc^{(m)})_{m\in\N}$, such that the corresponding subsequence of solutions $(u^{(m)})_{m\in\N}$ converges to $\bar u$ in $L^2(\O)$ and $(p^{(m)} - \half (u^{(m)})^2)_{m\in\N}$ weakly converges to $\bar p$ in $L^2(\O)$.
\end{theorem}

\medskip
\begin{proof}
Since the same estimates as in the linear case are available in the steady nonlinear case, the proof of proposition \ref{cvgce} holds for all the terms of \refe{nstocontfss} which are present in \refe{stocontf}.
We only have to prove that for a given $\varphi\in (C^\infty_c(\O))^d$, as $m\tends +\infty$:
\[
\terml{ss1}^{(m)} = \int_\O u^{(m)}(x)^2 \div_{\disc^{(m)}}(P_{\disc^{(m)}}\varphi)(x)\d x
\quad \mbox{ tends to } \quad
\int_\O \bar u(x)^2 \div\varphi(x)\d x
\]
and 
\[
\terml{ss2}^{(m)} = b_\disc(u^{(m)},u^{(m)},P_{\disc^{(m)}}\varphi)
\quad \mbox{ tends to } \quad
b(\bar u,\bar u,\varphi).
\]

\medskip
Thanks to the convergence in $L^2(\O)$ of $(u^{(m)})_{m\in\N}$ to $\bar u$ and to the
discrete Sobolev inequalities $\Vert v\Vert_{L^q(\O)} \le \ctel{sobdis} \Vert v\Vert_{\disc^{(m)}}$ 
for all $v\in H_{\disc^{(m)}}(\O)$ and all $q\le 6$ (see \cite[p. 790]{book}), 
we get using \refe{estimu} the convergence in $L^2(\O)$ of $((u^{(m)})^2)_{m\in\N}$ to $\bar u^2$.
We now remark that  for $i=1,\ldots,d$, the sequence $(P_{\disc^{(m)}}\varphi^{(i)})_{m\in\N}$ satisfies the hypotheses 
of Proposition \ref{cpct}. Hence, $\grad_{\disc^{(m)}} P_{\disc^{(m)}} \varphi^{(i)} $ weakly converges 
to $\grad \varphi^{(i)}$ in $L^2(\O)^d$. One has  $\div_\disc u = \sum_{i=1}^d
  \grad_{\disc}^{(i)} u^{(i)}$ for all 
$u \in (H_\disc(\O))^d$ such that $u_K = 0$ if $\edgescv \cap \edgesext \neq \emptyset.$
Hence $\div_{\disc^{(m)}}(P_{\disc^{(m)}}\varphi)$ weakly converges to $\div \varphi$ in $L^2(\O)$,
thus providing the limit of 
 $\termr{ss1}^{(m)}$.

\medskip
Thanks to \refe{deftridcbis}, setting for  simplicity  $\disc = \disc^{(m)}$,
we have:
\[
b_\disc(u,u,P_{\disc}\varphi)= \terml{ss3}^{(m)}-\terml{ss4}^{(m)}
\]
with:
\[
\ba
\termr{ss3}^{(m)}
& =& \dsp
\sum_{K\in\mesh} \sum_{L\in\NN_K} (A_{KL}\cdot u_K)((u_L - u_K) \cdot \varphi(x_K))
\\[3ex] & =& \dsp
\sumk \sumi \int_\O u^{(i)}(x) \grad_\disc^{(i)}(u^{(k)})(x) P_{\disc}\varphi^{(k)}(x)\d x
\\[4ex]
\termr{ss4}^{(m)}
& =& \dsp
\half \sum_{K|L\in \edgesint} (A_{KL}\cdot (u_L - u_K))((u_L - u_K) \cdot (\varphi(x_K) - \varphi(x_L)))
\ea
\]
Thanks to the convergence in $L^2(\O)$ of $(u^{(m)}P_{\disc^{(m)}}\varphi)_{m\in\N}$ to
$\bar u\varphi$, we get from proposition \ref{cpct} that:
\[
\lim_{m\to\infty} \termr{ss3}^{(m)} = \sumk \sumi \int_\O \bar u^{(i)}(x) 
\partial_i \bar u^{(k)}(x) \bar \varphi^{(k)}(x) \d x
= b(\bar u,\bar u,\varphi).
\]
We have:
\[
\ba
\termr{ss4}^{(m)} & = & \dsp \frac 1 4
\sum_{K|L\in \edgesint}d_{KL} (\tkl\n_{KL}\cdot (u_L - u_K))
((u_L - u_K) \cdot (\varphi(x_K) - \varphi(x_L)))
\ea
\]
and therefore, since $|\varphi(x_K) - \varphi(x_L)| \le d_{KL} C_\varphi \size(\disc)$
where $C_\varphi$ is a bound of $\grad\varphi$ in $L^\infty(\O)^d$, and since $d_{KL} \le 2\size(\disc)$, the following estimate holds:
\[
\ba
|\termr{ss4}^{(m)}|
& \le & \dsp 4 \size(\disc)^2 C_\varphi \Vert u\Vert_{\disc}^2.
\ea
\]
Therefore, \refe{estimu} yields:
\[
\lim_{m\to\infty} \termr{ss4}^{(m)} = 0,
\]
which concludes the proof of convergence.
\end{proof}

% -------------------------------------------------------------------------------------------

\subsection{The transient case}\label{trnl}

We now turn  to the study of the finite volume scheme for the transient Navier-Stokes equations, 
the weak formulation of which is given in \refe{weaksolt}.

We first give the definition of an admissible discretization for a space-time domain.

\begin{definition}[Admissible discretization, transient case]\label{adisct} 
Let $\O$ be an open bounded  polygonal (polyhedral if $d=3$) subset of $\R^d$, and
$\dr \O =  \overline{\O}\setminus\O$ its boundary, and let $T>0$.
An admissible finite volume discretization of $\O\times(0,T)$, denoted by $\disc$, is given
by $\disc=(\mesh,\edges,\centers,N)$, where $(\mesh,\edges,\centers)$ is an admissible
discretization of $\O$ in the sense of definition \ref{adisc} and $N\in\N_\star$ is given.
We then define $\dt = T/N$, and we denote by $\size(\disc) = \max(\size(\mesh,\edges,\centers),\dt)$ and $\regul(\disc) = \regul(\mesh,\edges,\centers)$.
\end{definition}

Under hypotheses \refe{hypomegat}-\refe{hypfgt}, let $\disc$ be an admissible discretization 
of $\O\times(0,T)$ in the sense of definition \ref{adisct} and let $\lambda\in(0,+\infty)$  
and $\alpha\in(0,2)$ be given.
We write a Crank--Nicholson scheme for the time discretization, and follow
 the  nonlinear steady--state case for the space discretization;
  the  finite volume scheme for the approximation of the 
solution \refe{nstocontt}--\refe{nstoconti} is then:
\be
\ba
u_0\in H_\disc(\O)^d, \\
\dsp u_{0,K} = \frac 1 {\mK} \int_K u_{\rm ini}(x) \d x,
\ \forall  K \in \mesh, 
\ea\label{schvfini}
\ee
and, again using Bernoulli's pressure $p + \half u^2$ instead of $p$, again denoted by $p$,
\be
\ba
(u_{n+1},p_{n+\half})\in (H_\disc(\O))^d\times H_\disc(\O), \\
\ \int_\O 
p_{n+\half}(x)\d x = 0,\
u_{n+\half} = \half(u_{n+1}+u_n),\\
\dsp \int_\O (u_{n+1}(x)-u_{n}(x))\cdot v(x) \d x + 
\dsp \nu \dt [u_{n+\half},v]_{\disc}  \\
\dsp - \dt \int_\O p_{n+\half}(x) \div_{\disc}(v)(x) \d x +
 \frac {\dt} 2 \int_\O u_{n+\half}(x)^2 \div_{\disc}(v)(x)\d x
\\ \dsp + \dt b_\disc(u_{n+\half},u_{n+\half},v)
 =   \dsp \int_{n\dt}^{(n+1)\dt}\int_\O f(x,t) \cdot v(x)\d x\d t, \\
\dsp \int_\O \div_\disc(u_{n+\half})(x) q(x) \d x  = - \lambda \ 
\size(\disc)^\alpha\   
\langle p_{n+\half},q\rangle_{\disc},  
\\ 
\dsp \ \forall  v \in H_\disc(\O)^d, 
\forall q\in H_\disc(\O),\ \forall n\in\N.
\ea\label{schvft}
\ee
In \refe{schvft}, we consider the approximation of $b_\disc$ given by \refe{deftridc}.
We then define the set $H_{\disc}(\O\times(0,T))$ of piecewise constant functions in each $K\times (n\dt,(n+1)\dt)$, $K\in\mesh$, $n\in\N$, and we define $(u,p)\in H_{\disc}(\O\times(0,T))$ by
\be
u(x,t) = u_{n+\half}(x),\ \hbox{ and }p(x,t) = p_{n+\half}(x),\ 
\hbox{ for a.e. }(x,t)\in\O\times(n\dt,(n+1)\dt),\ \forall n\in\N.
\label{solap}
\ee

\begin{remark}[Time discretization]
It is wellknown that the Crank--Nicholson discretization is implicit. 
If we use the $\theta$ scheme: $u_{n+\half} = \theta u_{n+1}+ (1-\theta) u_n$,
with $\theta\in[\half,1]$, 
the convergence proof which follows  applies  with  a few minor changes.
Variable time steps may also be considered.\end{remark}

Let us now prove the existence of at least one solution to scheme \refe{schvfini}-\refe{solap}.

\begin{proposition}[Existence of a discrete solution]
\label{exun} 
Under hypotheses \refe{hypomegat}-\refe{hypfgt}, let $\disc$ be an admissible
discretization of $\O\times(0,T)$ in the sense of Definition \ref{adisct}. 
Let $\lambda\in(0,+\infty)$  and $\alpha\in(0,2)$ be given.
Then there exists at least one $(u,p)\in (H_\disc(\O\times(0,T)))^d\times 
H_\disc(\O\times(0,T))$, solution to \refe{schvfini}-\refe{solap}.
\end{proposition}

\begin{proof}
We remark that, for a given $n=0,\ldots,N-1$, taking as unknown $u_{n+\half}$, 
and noting that $u_{n+1} = 2 u_{n+\half} - u_n$, 
Scheme \refe{schvft} is under the same form as scheme \refe{schvfnlss}, with $\eta = \frac {2} {\dt}$ 
and with a term in $u_n$ included in the
right hand side.
Therefore the existence of at least one solution follows from proposition \ref{exunss}.
\end{proof}

We then have the following estimate.

\begin{proposition}[Discrete $L^2(0,T;H^1_0(\O))$ estimate on velocities]\label{estl1l2t}
Under hypotheses \refe{hypomegat}-\refe{hypfgt}, let $\disc$ be an admissible discretization of 
$\O\times(0,T)$ in the sense of definition \ref{adisct}.
Let $\lambda\in(0,+\infty)$ and  $\alpha\in(0,2)$.
Let $(u,p)\in (H_\disc(\O\times(0,T)))^d\times H_\disc(\O\times(0,T))$, be a solution to 
\refe{schvfini}-\refe{solap}.
Then there exists $\ctel{0}>0$, only depending on $d$, $\O$, $\nu$, $u_0$, $f$, $T$ such 
that the following inequalities hold
\be
\Vert u\Vert_{L^\infty(0,T;L^2(\O)^d)} \le \cter{0},
\label{estimutinfv}\ee
\be
\Vert u\Vert_{L^2(0,T;H_\disc(\O)^d)} \le \cter{0},
\label{estimut}
\ee
and
\be
 \lambda \ \size(\disc)^\alpha \sum_{n=0}^{N-1} \dt
\vert  p_{n+\half}\vert_{\disc}^2 
= \lambda \ \size(\disc)^\alpha \int_0^T \vert 
 p(\cdot,t)\vert_\disc^2 \d t  
\le \cter{0}.
\label{estimpt}
\ee
\end{proposition}

\begin{proof}
Let $p=1,\ldots,N$.
We get, setting $v=u_{n+\half}$ in the first equation of \refe{schvft}, summing on $K\in\mesh$ and $n=0,\ldots,p-1$ in the first equation of \refe{schvft} and using property \refe{superb},
\[
\ba
\dsp \half \sum_{n=0}^{p-1}\int_\O (u_{n+1}(x)^2-u_{n}(x)^2) \d x + 
 \dsp \nu \sum_{n=0}^{p-1} \dt [u_{n+\half}, u_{n+\half}]_{\disc}
- \\ \dsp 
\sum_{n=0}^{p-1} \dt\int_\O p_{n+\half}(x) \div_{\disc}(u_{n+\half})(x) d 
x = \
 \sum_{n=0}^{p-1} \int_{n\dt}^{(n+1)\dt}\int_\O  f(x,t)\cdot 
u_{n+\half}(x)\d x\d t,
\ea
\]
This leads, setting $q = p_{n+\half}$  in the second equation of \refe{schvft}, to
\be
\ba
\dsp \half  \int_\O (u_{p}(x)^2-u_{0}(x)^2) \d x +
 \dsp \nu \sum_{n=0}^{p-1} \dt  [u_{n+\half}, u_{n+\half}]_{\disc}
+ \\ \dsp\lambda \ \size(\disc)^\alpha \sum_{n=0}^{p-1} \dt
\vert  p_{n+\half}\vert_{\disc}^2 = 
\dsp \int_{0}^{p\dt}\int_\O  f(x,t)\cdot u(x,t)\d x\d t.
\ea\label{ineqen}
\ee
Setting $p=N$ in  \refe{ineqen} gives \refe{estimut} and \refe{estimpt}.
The discrete Poincar\'e inequality \refe{poindis} and the inequality $\Vert u_0\Vert_{L^2(\O)^d} \le \Vert u_{\rm ini}\Vert_{L^2(\O)^d}$ give
\[
\ba
\dsp \Vert u_{p}\Vert_{L^2(\O)^d}^2 \le   
\frac {\diam(\O)^2} {2\nu}\Vert f\Vert_{L^2(\O\times(0,T))^d}^2 + \Vert 
u_{\rm ini}\Vert_{L^2(\O)^d}^2,\
\forall p = 1,\ldots,N,
\ea
\]
which proves \refe{estimutinfv}, since $\Vert u_{n+\half}\Vert_{L^2(\O)^d} \le
\half(\Vert u_{n}\Vert_{L^2(\O)^d}+\Vert u_{n+1}\Vert_{L^2(\O)^d})$ for all $n=0,\ldots,N-1$. \end{proof}

We then have the following estimates on translations.

\begin{proposition}[Space and time translate estimates]\label{transtime} 
Under hypotheses \refe{hypomegat}-\refe{hypfgt}, let $\disc$ be an admissible discretization of $\O\times(0,T)$ in the sense of definition \ref{adisct}.
Let $\lambda\in(0,+\infty)$, $\alpha\in(0,2)$  and $\theta >0$, such that $\regul(\disc) \ge \theta$.
Let $(u,p)\in (H_\disc(\O\times(0,T)))^d\times H_\disc(\O\times(0,T))$, be a solution to \refe{schvfini}-\refe{solap}.
We denote by $u$ the prolongment in $\R^d\times\R$ of $u$ by $0$ outside of $\O\times(0,T)$.
Then there exists $\ctel{1}>0$ and $\ctel{1t}>0$, only depending on $d$, $\O$, $\nu$, $\lambda$, $\alpha$, 
$u_0$, $f$, $\theta$ and $T$ such that the following inequalities hold:
\be
\Vert u(\cdot+\xi,\cdot)-u\Vert_{L^2(\R^d\times\R)}^2 \le \cter{1} 
|\xi|(|\xi| + 4\size(\mesh)),\
\forall \xi\in\R^d,
\label{estimtrsp}
\ee
and
\be
\Vert u(\cdot,\cdot+\tau)-u\Vert_{L^1(\R;L^2(\R^d))} \le
\cter{1t}|\tau|^{1/2},\
\forall \tau\in\R.
\label{estimtrti}
\ee
\end{proposition}

\begin{proof}
In the following proof, we denote by $C_i$, where $i$ is an integer, various positive real numbers
 which can only depend on $d$, $\O$, $\nu$, $\lambda$, $\alpha$, $u_0$, $f$, $\theta$  and $T$.
Inequality \refe{estimtrsp} is obtained from \refe{estimut} (see \cite{book}).
Let us prove \refe{estimtrti}.
Let $\tau\in(0,T)$ be given.
We define the following norms on $(H_\disc(\O))^d$, by:
\be
\ba
 \forall \  w\in (H_\disc(\O))^d, \\
 \Vert w \Vert_{\disc,\lambda}^2 =  \Vert w\Vert_{\disc}^2  +\\
  \dsp \frac 1 {\lambda\size(\disc)^\alpha}
\left(\sup\left\{\int_\O  \div_\disc(w)(x) q(x) \d x,\ q\in H_\disc(\O),
\vert   q\vert_\disc = 1\right\}\right)^2
\ea\label{normeel}
\ee
and
\be
\ba
 \forall \  w\in (H_\disc(\O))^d, \\
 \Vert w \Vert_{\star,\disc,\lambda} = \sup\left\{\int_\O  w(x)\cdot v(x) d
x, v\in (H_\disc(\O))^d,
\Vert v \Vert_{\disc,\lambda} = 1\right\}.
\label{normemel}
\ea
\ee
We then have, for a.e. $t\in(0,T)$,
\[
\Vert u(\cdot,t+\tau)-u(\cdot,t) \Vert_{L^2(\O)^d}^2 \le 
\Vert u(\cdot,t+\tau)-u(\cdot,t) \Vert_{\disc,\lambda}
\Vert u(\cdot,t+\tau)-u(\cdot,t) \Vert_{\star,\disc,\lambda},
\]
and therefore, thanks to the Young formula, 
\be\ba
\Vert u(\cdot,t+\tau)-u(\cdot,t) \Vert_{L^2(\O)^d} \le &
\frac {\sqrt{\tau}} 2 \Vert u(\cdot,t+\tau)-u(\cdot,t) 
\Vert_{\disc,\lambda}\\&
+ \frac 1 {2\sqrt{\tau}} \Vert u(\cdot,t+\tau)-u(\cdot,t) 
\Vert_{\star,\disc,\lambda}.
\label{teps}
\ea
\ee
We get, from \refe{schvft}, for all $q\in H_\disc(\O)$ and for a.e. $t\in(0,T)$,
\[
\int_\O  \div_\disc(u(\cdot,t))(x) q(x) \d x  = - \lambda \ \size(\disc)^\alpha\   
\langle p(\cdot,t),q\rangle_{\disc},
\]
which proves, using \refe{normeel}, that
\[
\Vert u(\cdot,t) \Vert_{\disc,\lambda}^2 \le \Vert u(\cdot,t)\Vert_{\disc}^2 + 
\lambda \ \size(\disc)^\alpha\vert p(\cdot,t)\vert_{\disc}^2.
\]
Using the Cauchy-Schwarz inequality, we have that:
\[
\left(\int_0^{T-\tau} \Vert u(\cdot,t+\tau)-u(\cdot,t) \Vert_{\disc,\lambda} \d 
t\right)^2 \le 
4 T \int_0^T \Vert u(\cdot,t) \Vert_{\disc,\lambda}^2 \d t,
\]
and therefore, using \refe{estimut}, \refe{estimpt},
\be
\int_0^{T-\tau} \Vert u(\cdot,t+\tau)-u(\cdot,t) \Vert_{\disc,\lambda} \d t 
\le \ctel{bdlam}.
\label{bordlambd}
\ee
We now study $\Vert u(\cdot,t+\tau)-u(\cdot,t) \Vert_{\star,\disc,\lambda}$.
We can write, for a.e. $t\in(0,T-\tau)$ and $x\in\O$,
\[
u(x,t+\tau)-u(x,t) =  
\half \sum_{n=0}^{N-1} (\chi_n(t,\tau)+ \chi_{n+1}(t,\tau)) (u_{n+1}(x) - 
u_{n}(x)), 
\]
where, for all $n\in\N$ and $t\in(0,T)$, $\chi_n(t,\tau) = 1$  if  $n\dt \in [t,t+\tau[$, and $\chi_n(t, \tau) =0$ otherwise.
This implies
\be
\ba\Vert u(\cdot,t+\tau)-u(\cdot,t) \Vert_{\star,\disc,\lambda} \le \\
\half \sum_{n=0}^{N-1} (\chi_n(t,\tau)+ \chi_{n+1}(t,\tau)) \Vert u_{n+1} 
- u_{n}\Vert_{\star,\disc,\lambda}. 
\ea\label{ght}
\ee
Let us then obtain a bound for $ \Vert u_{n+1} 
- u_{n}\Vert_{\star,\disc,\lambda}$.
Using the scheme \refe{schvft}, we get that, for all $v \in (H_\disc(\O))^d$,
\be
\begin{array}{ll}
\dsp \int_\O (u_{n+1}(x)-u_{n}(x))\cdot v(x) \d x & =    \dsp \int_{n\dt}^{(n+1)\dt}\int_\O 
f(x,t) \cdot v(x)\d x\d t  \\ \dsp
& -\dsp \nu \dt [u_{n+\half},v]_{\disc}  +     \dsp
\dt \int_\O p_{n+\half}(x)\div_{\disc}(v)(x) \d x  \\  & - \dsp
\frac {\dt} 2 \int_\O  u_{n+\half}^2 \div_{\disc}(v)(x) \d x -
\dt b_\disc(u_{n+\half},u_{n+\half},v) .
\end{array}
\label{coucou}\ee
Using the definition of $\div_{\disc}$, the fact that
$\sum_{\edge \in \edgescv} \medge \n_{K,\edge} = 0$,  and the
Cauchy-Schwarz inequality, there exists $\ctel{translat}$ 
such that:
\[
\int_\O u_{n+\half}^2(x) \div_{\disc}(v)(x) \d x \le \cter{translat} \Vert 
u_{n+\half}^2\Vert_{L^2(\O)} \ \Vert v \Vert_\disc.
\]
The discrete Sobolev inequality \refe{sob1} leads to
\[
\Vert u_{n+\half}^2\Vert_{L^2(\O)} \le \sumi \Vert  
(u_{n+\half}^{(i)})^2\Vert_{L^2(\O)} 
= \sumi \Vert u_{n+\half}^{(i)}\Vert_{L^4(\O)}^2  \le \ctel{ttyu} 
\Vert u_{n+\half}\Vert_\disc^{2}
\]
We take $\Vert v \Vert_{\disc,\lambda} = 1$ and note that, from Definition \refe{normeel},
we obtain that  $\Vert v \Vert_{\disc} \le 1$, and  that $\int_\O p_{n+\half}(x)\div_{\disc}(v)(x) \d x
 \le \left(\lambda \ \size(\disc)^\alpha\right)^{1/2}\   
\vert  p_{n+\half}\vert_{\disc}$.
We then pass  to the supremum in \refe{coucou}. 
Using the Cauchy-Schwarz inequality, the discrete Poincar\'e inequality, 
 and \refe{estib}, this yields:
\[
\begin{array}{ll}
\Vert u_{n+1} - u_{n}\Vert_{\star,\disc,\lambda} \le &
\dsp \sqrt{\dt}\diam(\O)  \Vert f \Vert_{L^2(\O\times(n\dt,(n+1)\dt))}\dt \\& +
\dsp  \dt \nu \Vert u_{n+\half} \Vert_{\disc}   
 + \left(\lambda \ \size(\disc)^\alpha\right)^{1/2}\   
\vert  p_{n+\half}\vert_{\disc} \\ & +\dsp
 \dt (\half  \cter{translat} \cter{ttyu} + \cter{bd})
\Vert u_{n+\half}\Vert_\disc^{2} .
\end{array}
\]
Summing the above equation for $n=0$ to $N-1$, applying the Cauchy-Schwarz inequality 
to all terms of the right hand side except the last, using \refe{estimut} 
and \refe{estimpt}, we get that there exists $\ctel{dvmu}$ such that
\[
\sum_{n=0}^{N-1} \Vert u_{n+1} - u_{n}\Vert_{\star,\disc,\lambda} \le 
\cter{dvmu}.
\]
Hence, noting that for all $n=0,\ldots,N$, $\int_0^{T-\tau} \chi_n(t,\tau) \d t\le \tau$,
we have:
\[
\half \int_0^{T-\tau} \sum_{n=0}^{N-1} (\chi_n(t,\tau)+ 
\chi_{n+1}(t,\tau)) 
\Vert u_{n+1} - u_{n}\Vert_{\star,\disc,\lambda} \d t 
\le  \cter{dvmu}\tau,
\]
which proves, using \refe{ght},
\be
\int_0^{T-\tau} \Vert u(\cdot,t+\tau)-u(\cdot,t) 
\Vert_{\star,\disc,\lambda} \d t\le  \cter{dvmu}\tau.
\label{machin}
\ee
Thanks to \refe{teps}, \refe{bordlambd} and \refe{machin}, we obtain that
\[
\int_0^{T-\tau} \Vert u(\cdot,t+\tau)-u(\cdot,t) \Vert_{L^2(\O)^d} \d t 
\le\ctel{cft} \sqrt{\tau}.
\]
Using \refe{estimutinfv}, we have
\[
\int_{T-\tau}^T \Vert u(\cdot,t+\tau)-u(\cdot,t) \Vert_{L^2(\O)^d} \d t = 
\int_{T-\tau}^T \Vert -u(\cdot,t) \Vert_{L^2(\O)^d} \d t \le \cter{0} \tau 
\le \sqrt{\tau}\sqrt{T}\cter{0},
\]
and a similar inequality holds for $\int_{-\tau}^0 \Vert u(\cdot,t+\tau)-u(\cdot,t) \Vert_{L^2(\O)^d} \d t$.
This thus gives \refe{estimtrti}, for any $\tau\in(0,T)$.
The case $\tau\ge T$ is obtained again using \refe{estimutinfv}, and the case $\tau\le 0$ is obtained from $\tau\ge 0$ by the change of variable $s = t+\tau$.
This completes the proof of \refe{estimtrti}.
\end{proof}

\begin{theorem}[Convergence of the  scheme]\label{cvgcetime}
Under hypotheses \refe{hypomegat}-\refe{hypfgt}, let $\theta >0$ be given and 
let $(\disc^{(m)})_{m\in\N}$ be a sequence of admissible discretizations 
of $\O\times(0,T)$ in the sense of definition \ref{adisct}, such that $\regul(\disc^{(m)}) \ge \theta$
and  $\size(\disc^{(m)})$ 
tends to $0$ as $m\to\infty$.
Let $\lambda\in(0,+\infty)$  and $\alpha\in(0,2)$ be given.
Let, for all $m\in\N$, $(u^{(m)},p^{(m)})\in (H_{\disc^{(m)}}(\O\times(0,T)))^d\times
H_{\disc^{(m)}}(\O\times(0,T))$, be a solution to \refe{schvfini}-\refe{solap} 
with $\disc = \disc^{(m)}$.
Then there exists a subsequence of  $(\disc^{(m)})_{m\in\N}$, 
again denoted $(\disc^{(m)})_{m\in\N}$, such that the corresponding subsequence 
of solutions  $(u^{(m)})_{m\in\N}$ converges in $L^2(\O\times(0,T))$ to a weak solution
$\bar u$ of \refe{nstocontt}-\refe{nstoconti} in the sense of definition \ref{weaksolt}.
\end{theorem}

\begin{proof}
Let us assume the hypotheses of the theorem.
Using translate estimates \refe{estimtrsp} and \refe{estimtrti} in the 
space $L^1(\R^d\times\R)$, we can apply Kolmogorov's theorem.
We get that there exists $\bar u \in L^1(\O\times(0,T))$ and a subsequence 
of $(\disc^{(m)})_{m\in\N}$, again denoted $(\disc^{(m)})_{m\in\N}$, such 
that the corresponding  subsequence of solutions  $(u^{(m)})_{m\in\N}$ converges in $L^1(\O\times(0,T))$ to $\bar u$ as $m\to\infty$.
Using \refe{estimut}, we get $\Vert u^{(m)} \Vert_{L^2(0,T;H_{\disc_m}(\O))} \le \cter{0}$, for all $m\in\N$, which gives, using the discrete Sobolev inequalities, $\Vert u^{(m)} \Vert_{L^1(0,T;L^4(\O))} \le \ctel{vbn}$, for all $m\in\N$.
Using a classical result on spaces $L^p(0,T;L^q(\O))$, we get that $(u^{(m)})_{m\in\N}$
converges in $L^1(0,T;L^2(\O))$ to $\bar u$ as $m\to\infty$.
Thanks to \refe{estimutinfv}, we have $\Vert u^{(m)}\Vert_{L^\infty(0,T;L^2(\O)^d)} \le \cter{0}$, for all $m\in\N$.
The same result on spaces $L^p(0,T;L^q(\O))$ implies that $(u^{(m)})_{m\in\N}$ converges
in $L^2(0,T;L^2(\O))$ to $\bar u$ as $m\to\infty$.
We can therefore pass to the limit in \refe{estimtrsp}.
The resulting inequality implies $\bar u \in L^2(0,T;H^1_0(\O)^d)$ (see \cite{book}).
Passing to the limit in \refe{estimutinfv} leads to $\bar u \in L^\infty(0,T;L^2(\O)^d)$.

Let us now prove that $\bar u$ is a weak solution of \refe{nstocontt}-\refe{nstoconti} in the sense of definition \ref{weaksolt}.

Let $\varphi\in C_c^\infty(\O\times(-\infty,T))^d$ be given, with $\div\varphi(x,t) = 0$ for all $(x,t)\in\O\times(-\infty,T)$.
Let $\disc^{(m)}$ be a given admissible discretization extracted from the considered subsequence.
Omitting some of the indices $m$ for the simplicity of notation, 
we then set $v = P_\disc\varphi(\cdot,n\dt)$ in \refe{schvft},
 and we sum for $n=0,\ldots,N-1$.
We thus get
\be
\termr{t1}^{(m)}+\termr{t2}^{(m)}+\termr{t3}^{(m)}+\termr{t3b}^{(m)}+\termr{t4}^{(m)} = \termr{t5}^{(m)},
\label{sumtt}
\ee
with
\[
\terml{t1}^{(m)} = \sum_{n=0}^{N-1} \int_\O (u_{n+1}(x)-u_{n}(x)) \cdot 
P_\disc\varphi(x,n\dt) \d x,
\] 
\[
\terml{t2}^{(m)} = \sum_{n=0}^{N-1} \dt [u_{n+\half}, P_\disc\varphi(\cdot,n\dt)]_{\disc},
\] 
\[
\terml{t3}^{(m)} = -\sum_{n=0}^{N-1} \dt \int_\O p_{n+\half}(x) 
\div_{\disc}(P_\disc\varphi(\cdot,n\dt))(x) \d x,
\] 
\[
\terml{t3b}^{(m)} = \half \sum_{n=0}^{N-1} \dt \int_\O u_{n+\half}(x)^2 
\div_{\disc}(P_\disc\varphi(\cdot,n\dt))(x) \d x,
\] 
\[
\terml{t4}^{(m)} = \sum_{n=0}^{N-1} \dt 
b_\disc(u_{n+\half},u_{n+\half},P_\disc\varphi(\cdot,n\dt)),
\]
and 
\[
\terml{t5}^{(m)} = \sum_{n=0}^{N-1} \int_{n\dt}^{(n+1)\dt}\int_\O  f(x,t)\cdot 
P_\disc\varphi(x,n\dt)\d x\d t.
\] 
In the following, we denote by $C_i$ various positive reals which 
can only depend on $d$, $\O$, $T$, $u_{\rm ini}$, $f$, $\nu$, $\theta$ and $\lambda$.
We first start with the study of $\termr{t2}$.
We classically have (see \cite{book})
\be
\lim_{m\to\infty} \termr{t2}^{(m)} = \int_0^T \int_\O \grad \bar 
u(x,t):\grad \varphi(x,t) \d x \d t.
\label{cvt2}
\ee
The proof that
\be
\lim_{m\to\infty} \termr{t3}^{(m)} = 0
\label{cvt3}
\ee
is a consequence of \refe{estimpt} and of a direct adaptation of Proposition
\ref{weakconvgrad} to time-dependent functions.
Let us now prove that
\be
\lim_{m\to\infty} \termr{t3b}^{(m)} = 0.
\label{cvt3b}
\ee
Since $(u^{(m)})^2$ tend to $\bar u^2$ as $m\to\infty$ in $L^1(\O\times(0,T))$, 
the same argument as in the steady state case (see proof of theorem \ref{cvgcenlss}) 
provides \refe{cvt3b}.

We now turn to the study of $\termr{t4}$.
Following the proof of proposition \ref{cvgcenlss}, the proof that
\be
\lim_{m\to\infty} \termr{t4}^{(m)} = \int_0^T  b(\bar u(\cdot,t),\bar 
u(\cdot,t),\varphi(\cdot,t)) \d t.
\label{cvt4}
\ee
is a direct consequence of the convergence of $u$ to $\bar u$ in $L^2(\O\times(0,T))$ and proposition \ref{weakconvgrad}.
The study of $\termr{t5}$ is classical, and we have
\be
\lim_{m\to\infty} \termr{t5}^{(m)} = \int_0^T \int_\O f(x,t) 
\cdot\varphi(x,t) \d x \d t.
\label{cvt5}
\ee
Let us now prove that
\be
\lim_{m\to\infty} \termr{t1}^{(m)} = - \int_0^T  \int_\O \bar 
u(x,t)\partial_t \varphi(x,t) \d x \d t - 
\int_\O  u_{\rm ini}(x)\varphi(x,0) \d x.
\label{cvt1}
\ee
Indeed, we have
\[
\termr{t1}^{(m)} = - \int_\O u_{0}(x) \cdot P_{\disc}\varphi(x,0) \d x
- \terml{t1deplus}^{(m)} - \frac 1 2 \terml{t1b}^{(m)}.
\]
with 
\[
\termr{t1deplus}^{(m)} =  
  \sum_{n=0}^{N-1} \int_\O u_{n+\frac  1 2}(x)\cdot (P_{\disc}\varphi(x,(n+1)\dt) - 
P_{\disc}\varphi(x,n\dt)) \d x.
\]
and 
\[
\termr{t1b}^{(m)} = \sum_{n=0}^{N-1} \int_\O (u_{n+1}(x) - u_n(x)) \cdot 
(P_{\disc}\varphi(x,(n+1)\dt) - P_{\disc}\varphi(x,n\dt)) \d x
\]
We classically have
\[
\lim_{m\to\infty} \int_\O u_{0}(x) \cdot P_{\disc}\varphi(x,0) \d x = \int_\O  
u_{\rm ini}(x)\varphi(x,0) \d x.
\]
We also easily have, thanks to the convergence properties of $u^{(m)}$, that
\[
\lim_{m\to\infty}  \termr{t1deplus}^{(m)}  =
\int_0^T  \int_\O \bar u(x,t)\partial_t\varphi(x,t) \d x \d t.
\]
Let us prove that the term $
\termr{t1b}^{(m)}$  
tends to 0 as $m\to\infty$. We have 
 $\termr{t1b}^{(m)}= \terml{t1bp}^{(m)} - \termr{t1}^{(m)}$, with
\[
\termr{t1bp}^{(m)} = \sum_{n=0}^{N-1} \int_\O (u_{n+1}(x) - u_n(x)) \cdot 
P_{\disc}\varphi(x,(n+1)\dt) \d x.
\]
Thanks to the limits given by \refe{cvt2}, \refe{cvt3}, \refe{cvt3b}, 
\refe{cvt4} and \refe{cvt5}, and thanks to \refe{sumtt}, we obtain that
$\dsp \lim_{m\to\infty}\termr{t1}^{(m)} = \terml{lim}$, with
\[\ba
\termr{lim} = &
-\dsp\nu \sumi \int_0^T\int_\O \nabla u^{(i)}(x,t)\cdot\nabla 
\varphi^{(i)}(x,t)\d x \d t -
\int_0^T b(u(\cdot,t),u(\cdot,t),\varphi(\cdot,t)) \d t + \\ &
\dsp  \int_0^T\int_\O f(x)  \cdot\varphi(x,t) \d x \d t.
\ea
\]
Since \refe{cvt2}, \refe{cvt3}, \refe{cvt3b}, \refe{cvt4} and \refe{cvt5} are available as well, replacing $P_{\disc}\varphi(\cdot,n\dt)$ by $P_{\disc}\varphi(\cdot,(n+1)\dt)$
in $\termr{t2}$, $\termr{t3}$, $\termr{t3b}$, $\termr{t4}$ and  $\termr{t5}$, we also get using \refe{schvft} with $v=P_{\disc}\varphi(\cdot,(n+1)\dt)$ that
$
\dsp \lim_{m\to\infty}\termr{t1bp}^{(m)} = \termr{lim}. 
$
Thus we get that  $\lim_{m\to\infty}\termr{t1b}^{(m)} = 0$, which concludes the proof of \refe{cvt1}.
Thanks to  \refe{sumtt}, \refe{cvt1}, \refe{cvt2}, \refe{cvt3}, \refe{cvt3b}, 
\refe{cvt4} and \refe{cvt5}, we thus obtain \refe{nstocontft}, provided that we can prove
\[
\div\bar u(x,t) = 0,\ \hbox{ for a.e. }(x,t)\in\O\times(0,T).
\]
This last relation can be shown, following the proof of \refe{vandiv}.
This completes the proof of the above theorem.
\end{proof}

\begin{remark}
Using the above proof of convergence, we get the energy inequality for $d=2$ or $3$ from inequality \refe{ineqen}, since we have the property
\[
\int_0^T\int_\O (\grad u^{(i)}(x,t))^2 \d x \d t \le \liminf_{m\to\infty} 
\sum_{n=0}^{N^{(m)}-1} \dt  [u_{n+\half}^{(m,i)}, 
u_{n+\half}^{(m,i)}]_{\disc^{(m)}}
\] 
\end{remark}

\section{Numerical results}\label{secnum}
An industrial implementation of a colocated finite volume scheme
may be found in \cite{neptune} for instance, 
where complex applications are considered.
Focusing in this paper on properties of convergence and error estimates,   
some simple numerical experiments are described here to observe 
the convergence rate of Schemes \refe{schvf} and \refe{schvfini}-\refe{schvft}
with respect to 
the space and time discretizations.  To that purpose, we use a prototype code
where the nonlinear equations  are solved by an underrelaxed Newton  method, and
 the linear systems by a direct band Gaussian elimination solver. This code  
handles Stokes or Navier-Stokes
problems with various boundary conditions, 
using non uniform rectangular or triangular meshes on general 2D polygonal domains.

\medskip

The linear Stokes equations are first considered in the case
$d=2$,  $\O = (0,1)\times(0,1)$, $\nu = 1$,  and
$f$ is taken to satisfy \refe{stocont} with a solution equal to
\[\ba
\bar u^{(1)}(x^{(1)},x^{(2)}) &= -  \dr^{(2)} \Psi( x^{(1)},x^{(2)})\\
\bar u^{(2)}(x^{(1)},x^{(2)}) &=  \dr^{(1)} \Psi( x^{(1)},x^{(2)}) \\
\bar p(x^{(1)},x^{(2)}) &= 100 \  \left((x^{(1)})^2 + (x^{(2)})^2\right),
\ea\]
denoting by $\Psi( x^{(1)},x^{(2)}) = 1000 \ [x^{(1)}(1-x^{(1)})x^{(2)}(1-x^{(2)})]^2$.
The approximate solution $(u,p)$ is computed with the scheme \refe{schvf}.
The observed numerical
order of convergence, considering the norms 
$\Vert u - P_\disc\bar u\Vert_{L^2(\O)}$
and $\Vert p - P_\disc\bar p\Vert_{L^2(\O)}$, is equal to $2$ for the velocity components, 
and to $1$ for the pressure in
the cases of non uniform rectangular and square meshes (from 400 to 6400 grid blocks).
Note that in these cases, 
there is apparently no need for
a significant positive value of the stabilization coefficient $\lambda$.  
The observed numerical 
order of convergence is similar in the case of triangular meshes (from 1400 to 5600 grid blocks), 
but values such as $\lambda = 10^{-4}$, $\alpha = 1$ have to be
used in order to avoid oscillations in the pressure field. This confirms that in the 
case of triangles, the approximate pressure space is too large to avoid 
stabilization. In fact, other tests were performed (e.g. the classical backward step) 
which show that stabilization is also needed in the case of rectangles when more 
severe problems are considered. Note that in industrial implementations, stabilization 
may be performed with other means, see \cite{fluent}, \cite{neptune}, (see also
\cite{boivin} in the triangular case). 

\medskip

We then proceed to a similar comparison in the case of transient nonlinear problems.
Considering a transient adaptation of the above
steady-state analytical solution, the continuous problem is 
then defined by zero initial and
boundary conditions, $T = 0.1$, and the function
$f$ is taken to satisfy \refe{nstocontt} with a solution equal to
\[\ba
\bar u^{(1)}(x^{(1)},x^{(2)},t) &= -t \ \dr^{(2)} \Psi( x^{(1)},x^{(2)})\\
\bar u^{(2)}(x^{(1)},x^{(2)},t) &= t \ \dr^{(1)} \Psi( x^{(1)},x^{(2)}) \\
\bar p(x^{(1)},x^{(2)},t) &= 100 \ t\ \left((x^{(1)})^2 + (x^{(2)})^2\right),
\ea\]
with the same function $\Psi$ as above.
We again observe an  order 2 of convergence of the approximate solution
at times $t = .05$ and $t = .1$,
when the space and the time discretizations are simultaneously modified
with the same ratio (from $\dt = 0.01$ to $\dt = 0.0025$ as the size of the mesh
is divided by 4). Similar observations are still valid for 
the classical Green-Taylor example.
%is given for $d=2$ and $\O = (0,1)\times(0,1)$, $T = 1/10$, by $f=0$, nonhomogeneous Dirichlet
%boundary conditions and
%\[\ba
%\bar u^{(1)}(x^{(1)},x^{(2)},t) &= -\cos(2\pi(x^{(1)}+\frac 1 4))\ \sin(2\pi(x^{(2)}+\frac 1 2))\ \exp(-8\pi^2 t
%\nu) \\
%\bar u^{(2)}(x^{(1)},x^{(2)},t) &= \sin(2\pi(x^{(1)}+\frac 1 4))\ \cos(2\pi(x^{(2)}+\frac 1 2))\ \exp(-8\pi^2 t
%\nu) \\
%\bar p(x^{(1)},x^{(2)},t) &= -\frac 1 4 \left( \cos(4\pi(x^{(1)}+\frac 1 4))+
%\cos(4\pi(x^{(2)}+\frac 1 2))\right)\ \exp(-16 \, \pi^2 t
%\nu).
%\ea\]

\section{Conclusions}\label{seconcrem}

The above numerical results show that the  theoretical 
error estimate which is proved in Section \ref{secfvslin} 
for the linear Stokes equations is non optimal; a sharper estimate is currently 
being written \cite{EHL} under more regularity assumptions on the mesh.  

The proof of convergence of the full  space-time discrete approximation of 
\refe{nstocontt}  given by \refe{schvft}  uses  estimates on the time 
translates, which were introduced in the $L^2(\Omega\times (0,T))$ framework
for the proof of convergence 
 of the finite volume method for degenerate parabolic equations \cite{slimane,book} 
 and used for several other cases, see e.g.\cite{convpardeg}. A major difficulty 
 which arises here is the handling on the nonlinear convective term, as in the
 continuous case, which leads us to establish an estimate on the time 
translates in  $L^1(0,T; L^2(\Omega))$. This new technique may be used for 
parabolic problems with other type of nonlinearities. 

We remarked that industrial codes use other types of stabilizations than the one 
used here. Further works will be devoted to the mathematical study of 
such stabilizations, for which, to our knowledge, 
no proof of convergence is known up to now. 

Finally, let us also mention undergoing work on a generalization of the scheme studied here
to  the  full transient Navier-Stokes equations 
including the energy balance,  under the 
Boussinesq approximation.

\end{document}